\documentclass[final,3p, authoryear]{elsarticle}

\usepackage{amssymb}
 \usepackage{amsthm}

\usepackage{longtable}
\usepackage{amsmath}
\usepackage{algorithm}
\usepackage{algorithmic}
\usepackage {ctable}

\usepackage{nomencl}
\makenomenclature
\usepackage{etoolbox}

\usepackage[colorlinks=true,breaklinks=true,bookmarks=true,urlcolor=blue,
citecolor=blue,linkcolor=blue,bookmarksopen=false,draft=false]{hyperref}

\newcommand{\T}{\mathcal{T}}
\newcommand{\Ccal}{\mathcal{C}}
\newcommand{\Ecal}{\mathcal{E}}
\newcommand{\Ncal}{\mathcal{N}}
\newcommand{\Sbb}{\mathbb{S}}
\newcommand{\Scal}{\mathcal{S}}
\newcommand{\Jbb}{\mathbb{J}}
\newcommand{\Ibb}{\mathbb{I}}

\newcommand{\y}{\boldsymbol{y}}
\newcommand{\bb}{\boldsymbol{b}}

\newcommand{\sto}{\mbox{s.~t.}}

\usepackage{color}

\xdefinecolor{auburnORANGEdark}{HTML}{dd550c}

\def\orange#1{{ #1}}

\newdefinition{definition}{Definition}
\newproof{prf}{Proof}
\newdefinition{remark}{Remark}

\newtheorem{theorem}{Theorem}

\newtheorem{lemma}[theorem]{Lemma}
\begin{document}

\begin{frontmatter}

%% Title, authors and addresses

%% use the tnoteref command within \title for footnotes;
%% use the tnotetext command for theassociated footnote;
%% use the fnref command within \author or \address for footnotes;
%% use the fntext command for theassociated footnote;
%% use the corref command within \author for corresponding author footnotes;
%% use the cortext command for theassociated footnote;
%% use the ead command for the email address,
%% and the form \ead[url] for the home page:
%% \title{Title\tnoteref{label1}}
%% \tnotetext[label1]{}
%% \author{Name\corref{cor1}\fnref{label2}}
%% \ead{email address}
%% \ead[url]{home page}
%% \fntext[label2]{}
%% \cortext[cor1]{}
%% \address{Address\fnref{label3}}
%% \fntext[label3]{}

\title{A Graph Theoretic Approach to Non-Anticipativity Constraint Generation in Multistage Stochastic Programs with Incomplete Scenario Sets}

%% use optional labels to link authors explicitly to addresses:
 \author[AUCE]{Brianna Christian}
 \author[AUIE]{Alexander Vinel \corref{cor1}}
 \author[AUCE]{Zuo Zheng}
 \author[AUCE]{Selen Cremaschi}
\address[AUCE]{Department of Chemical Engineering, Auburn University}
\address[AUIE]{Department of Industrial and Systems Engineering, Auburn University}
\cortext[cor1]{alexander.vinel@auburn.edu}

%% \address[label2]{}

\begin{abstract}
%% Text of abstract
We propose an algorithm for generating a minimum-cardinality set of non-anticipativity constraints (NAC) for scenario-based multistage-stochastic programming (MSSP) problems with both endogenous and exogenous uncertainties which allow for  gradual realizations. 
Recently several authors have considered approaches to generate the minimum cardinality NAC set for MSSPs for various scenario set structures. However, these approaches have been limited to uncertain parameters where the realizations occur instantaneously or the full set of scenarios is required. The proposed algorithm, referred to as Sample Non-Anticipativity Constraint algorithm (SNAC) relaxes this requirement. We show that as long as the number of uncertain parameters and parameter values are kept constant, the algorithm scales polynomially in the number of scenarios. 

\end{abstract}

\begin{keyword}
%% keywords here, in the form: keyword \sep keyword
Stochastic programming \sep Multistage Stochastic Models \sep Endogenous Uncertainty \sep Non-Anticipativity Constraints 
%% PACS codes here, in the form: \PACS code \sep code

%% MSC codes here, in the form: \MSC code \sep code
%% or \MSC[2008] code \sep code (2000 is the default)

\end{keyword}

\end{frontmatter}

%% \linenumbers
\sloppy
%% main text
\section{Introduction}\label{intro} %%1.
%\cite{apap-2017}

Stochastic programming 	is a scenario-based approach for optimization under uncertainty \citep{birge-2011}.
%, and it has been applied to several problems including among others portfolio optimization \citep{Vigerske-2007, Solak-2010}, water resource management \citep{li-2008}, and energy production \citep{Dentcheva-1998, Fleten-2008}. 
%The approach is characterized by a set of discrete stages where decisions are made and uncertainty is realized, and scenarios which represent the potential outcomes of uncertainty. 
%Stochastic programs are classified based on the number of stages. 
Two-stage stochastic programs have a first stage where decisions are made and a second stage where uncertainty is realized and recourse actions can be taken. In contrast, multistage stochastic programs (MSSPs) have more than one stage where decisions are made and one or more periods where uncertainty is realized and recourse actions are taken. Note that the number of time periods in the planning horizon does not necessarily correspond to the number of stages in the problem.
%Figure 1 highlights the difference between two-stage and MSSPs. It is important to note that the number of time periods in the planning horizon do not necessarily correspond to the number of stages in the problem. \red{\emph{I don't think the figure here adds much}} 

%\begin{figure}
%	\FIGURE
%	{\includegraphics[width=0.5\textwidth]{two-stage-model.png}}
%{Comparison of (a) two-stage stochastic programs and (b) multi-stage stochastic programs \label{fig:two-stage-model}}
%{}
%\end{figure}

%The second component of a stochastic programming problem is a finite set of representative scenarios. 
The scenarios in a stochastic program are generated based on the outcomes of the uncertain parameters. A common approach for scenario generation in problems where uncertain parameters can be considered independent is to take the Cartesian product of the outcomes of the uncertain parameters, as described by  \citet{apap-2017}. 
%In situations where the distribution of the uncertain parameter is continuous, a discrete distribution is often used as an approximation. 
Uncertainty in stochastic programming can be classified based on how and when it is realized. Realization of \emph{exogenous}, or \emph{decision-independent}  uncertain parameters occurs at a known time independent of the values of the decision variables \citep{birge-2011}. In operations research, the market conditions such as demand are often assumed to be independent of any investment decisions and are treated as exogenous uncertain parameters. 
%The other type of uncertainty is endogenous uncertainty. 
With \emph{endogenous}, or \emph{decision-dependent} uncertainty, the timing of the realization or the value of the realized parameter is impacted by the value of the decision variables \citep{Jonsbraten-1998}. 
%This type of uncertainty can be further classified based on whether the values of the decision variables impact only the timing of the realization (i.e. endogenous uncertainty with exogenous realization) or whether both the timing and the value that is realized are impacted by the values of the decision variables. An example of the first type of endogenous uncertainty is encountered in the clinical trial planning problem \citep{colvin-2008}.  Uncertainty in whether a clinical trial is successful or not is only realized after the clinical trial has been completed, however whether the clinical trial is successful or not is not impacted by when the clinical trial is started. In the second type of endogenous uncertainty, the values of the decision variables impact the timing of the realization and/or the value of the realized uncertainty. In literature, this type of endogenous uncertainty can arise in facility protection problem, where, the likelihood of a facility failing to deliver goods or services after a disruptive event depends on the level of resources allocated as protection to that facility \citep{medal-2016}. 

Another way to classify uncertainty is based on how the uncertainty is realized. Uncertainty in stochastic programs can either be realized instantaneously or gradually along the planning horizon \citep{solak-2007}. In the case of exogenous uncertainty, whether it is realized instantaneously or gradually has little to no effect on the multistage stochastic programming formulation, which is not the case for endogenous parameters. A well-known example of this type of uncertainty realization can be seen in the planning and scheduling of pharmaceutical clinical trials \citep{colvin-2008}. In the problem, the goal is to find the clinical trial schedule which provides the highest expected net present value given uncertainty in the outcome of clinical trials. The authors model the clinical trial outcome uncertainty for each potential new product using a single uncertain parameter. Uncertainty is only realized if a clinical trial fails or if all clinical trials are successfully completed. This leads to time periods in the planning horizon where a product may have successfully completed one or more clinical trials (i.e., no failures) but has not successfully completed all of them. 
%In MSSPs, partial realizations of uncertain endogenous parameters complicate the definition of constraints sets which ensure that the underlying value of the uncertain parameter is not anticipated.

The formulation of an MSSP requires the equality of decision variables in different scenarios  to ensure that the realized values of the uncertain parameters are not anticipated. Depending on the type of uncertainty in the problem, this can be accomplished either implicitly or explicitly. When the timing of the realization of uncertainty is known (i.e., exogenous uncertainty), anticipation of the underlying value of the uncertain parameters can be prevented by defining a variable used to relate the decision variables in two different scenarios. Using the same decision variable in multiple scenarios ensures that the decision variable values are identical in scenarios where uncertainty has not been realized. Alternatively, a set of equality constraints, called non-anticipativity constraints (NACs), which set the values of the decision variables equal to each other at time periods before uncertainty is realized can be used. If the time period of differentiation is not known (i.e., endogenous uncertainty), a set of NACs must be used as there is no way to implicitly incorporate the realization of uncertainty.

In this paper, we introduce and describe in detail a graph theoretic approach that generates a minimum cardinality NAC set for MSSPs with an incomplete scenario set where the realization of uncertainty occurs gradually. 
\orange{
While in principle it  is possible to reformulate the problem  as an instantaneous realization model by introducing additional uncertain parameters (and consequently apply some of the existing methods), our approach has the benefit of operating directly on the natural original formulation. 
}

	The remainder of the paper is organized as follows. In Section \ref{sec:lit} we review existing relevant literature. Next, we present the notation and give the general NAC formulation for MSSPs with gradual realization of uncertainty in Section \ref{sec:form}.
	In Section \ref{sec:method}, we propose the algorithm, which generates a minimum cardinality NAC set. We also present an illustrative example for a small-scale problem (Section \ref{sec:comp-simple}) and report results of a computational case study (Section \ref{sec:comp}). Conclusions and future directions are given in Section \ref{sec:conclusions}.

%For the sake of completeness, in Section \ref{sec:background} we present the general standard formulations of MSSPs with endogenous (Section \ref{sec:en}) and with both endogenous and exogenous uncertainty (Section \ref{sec:enex}). 
%The formulation for the MSSP with just endogenous uncertainty is given in , and the extension to MSSPs with both endogenous and exogenous uncertainty is shown in . 

\subsection{Review of the Literature on NAC Generation}\label{sec:lit}

Multistage stochastic programming problems that include endogenous uncertain parameters is a relatively new research area and has received far less attention in literature compared to its counterpart under exogenous uncertainties \citep{apap-2017}. To the best of our knowledge, \cite{Jonsbraten-1998} were the first  to introduce endogenous uncertainty by extending size selection problem, which was defined originally by \cite{jorjani}. 
Since then there has been an increasing number of studies applying endogenous uncertainty in gas-field development problem  \citep{goel2004, Goel-2006}, in synthesis of process networks \citep{tarhan2008}, in open pit mining scheduling problem by \citep{boland2008}, in a pharmaceutical clinical trial planning problem \citep{colvin-2008}, in offshore oilfield development planning \citep{gupta2014}, in vehicle routing problem \citep{RN9}, in various portfolio optimization of smart grid technologies \citep{RN10}, in artificial lift infrastructure planning problem \citep{RN11}, and in a new technology investment planning (NTIP) problem by \citep{RN12}.

However, application of MSSPs is limited due to curse of dimensionality. Rapid  growth of the number of scenarios gives rise to a large number of NACs requiring  a large number of additional indicator variables and logical constraints. Recent literature in the area focuses on the solution approaches to address the space and time complexity of MSSP models, especially due to the challenges associated with the large number of NACs. 

Lagrangean relaxation approaches have received a significant attention. Here, NACs are removed or relaxed partially, which decomposes the MSSP into smaller problem of scenario groups. \cite{gupta2014} proposed a Lagrangean decomposition algorithm to obtain dual bounds for MSSP models under endogenous uncertainties. The authors decomposed scenarios into different scenario groups, only dualized the initial NACs linking two scenario groups, and kept initial and conditional NACs for the scenario pairs in the same scenario group. They applied this new Lagrangean decomposition algorithm to process-network-synthesis problems and reached 1\% optimality gap within 30 iterations. \cite{RN5} presented a branch-and-bound algorithm that generates upper-bounds using the solution of the Lagrangean dual problem with relaxed NACs. The lower bounds were generated heuristically based on the solution of the Lagrangean dual problem. The results suggested that the branch-and-bound algorithm achieved significantly better solutions and tighter optimality gaps than the heuristic presented in \cite{goel2004}. \cite{RN14} presented an approach that improves dual bounds for solving nonlinear MSSPs with decision-dependent uncertainty. The algorithm solves the Lagrangean relaxation of the dual problem where NACs have been removed to obtain the dual bound. The primal bound is generated by locating a feasible solution using a rolling-horizon approach. 

Another approach utilizes scenario decomposition or feasible solution space decompositions. \cite{goel2004} considered an approximation approach, which searched a sub-space of the feasible region to find a ``good'' solution. They used this approach for solving a gas-field development problem, where the size and initial deliverability of reserves were endogenous uncertainties realized immediately after the drilling decisions were made. The sub-space is obtained by removing the scenario dependency of investment decisions for uncertain fields, yielding a more constrained version of the original problem. Then, this constrained version is solved to optimality. The results revealed that the solutions obtained to the constrained formulation were significant improvements compared to the deterministic solutions. \cite{Solak-2010} developed a sample average approximation (SAA) based algorithm to generate candidate solutions for the R\&D project portfolio optimization problem, where return levels for projects were endogenous uncertain parameters. \cite{apap-2017} developed a sequential scenario decomposition (SSD) approach for MSSP models under both endogenous and exogenous uncertainties. The heuristic first solves endogenous sub-problems by selecting one scenario from each exogenous scenario group. Each sub-problem includes first-period NACs and endogenous NACs, but they do not include exogenous NACs. The approach solves a series of sub-problems with endogenous uncertainties to determine binary investment decisions, fixes the binary investment decisions to satisfy the first-time period endogenous and exogenous NACs, and solves the resulting model to obtain a feasible solution. \cite{apap-2017} applied SSD to the process synthesis problem \citep{tarhan2008, Goel-2006} and the offshore oilfield planning problem \citep{RN18}, and obtained high-quality feasible solutions.

In terms of branch and bound approaches, \cite{tarhan2008} proposed a solution strategy that used a duality based branch-and-bound algorithm for a MSSP model for solving process network synthesis problem, where the process yields are endogenous uncertain parameters with gradual realization as process operating decisions are made. The strategy was able to solve a 16 scenario example to within 3\% optimality. \cite{RN18} explored an offshore oil planning problem with endogenous uncertainty in the initial maximum oil flowrate, and considered a single oil field with non-linear reservoir model and a gradual realization of uncertainty. Using a duality-based branch-and-bound, the authors located solutions that were up to 22\% better than solutions obtained using an expected value approach. They, however, noted that the solution times for the MSSP were ``rather long''. \cite{RN19} developed a branch-and-bound algorithm combining a heuristic primal bounding approach, Knapsack problem based Decomposition Algorithm (KDA) \citep{RN20} and two dual-bounding approaches to solve large-scale MSSPs with endogenous uncertainty. The dual bounds were generated using a modified approach based on progressive hedging (PH), and by solving the individual scenario problems. The primal bounds were generated using the KDA. The branch-and-bound algorithm successfully addressed the space complexity of MSSP problems, but required considerable time to converge.

\cite{colvin-2008} introduced a MSSP model for planning clinical trials in pharmaceutical R\&D pipelines. The authors explored a rolling-horizon approximation approach, which yielded tight feasible solutions for large instances of the problem \citep{RN22}, where the planning horizon was divided into a finite number of subsets. A relaxed MSSP is generated for the first subset by removing all inequality NACs for the stages beyond the first subset. The solution of the relaxed MSSP is implemented for the first subset, and related uncertainty is realized. Then, the process is repeated for each subset until the end of the planning horizon is reached. The authors were able to successfully solve cases with more than 1000 scenarios. Another approach, a branch-and-cut algorithm \citep{RN23}, initially adds a percentage of NACs and then iteratively adds constraints based on violations. The authors concluded that the algorithm reduced the number of NACs that should be included in the problem formulation significantly and that this method would be advantageous for any problem where the majority of constraints are NACs. 

\cite{RN24} introduced a novel approach that uses decision rules for approximating MSSP problems under endogenous uncertainties. The approach employs piecewise constant and piecewise linear functions of uncertainties and can handle uncertain parameters with both discrete and continuous distributions. The authors applied the approach to solve oil infrastructure planning problem, and the approach yielded improved primal bounds with limited time budget.

In order to address the computational complexity caused by NACs and large size of scenarios, several authors considered approaches that remove the scenario and NAC structure to obtain good feasible solutions with limited computational effort. \cite{RN25} proposed a multi-step anticipatory algorithm, which uses a sample average approximation approach to generate a Markov Decision Process (MDP). The MDP is then solved, and the greedy solution is returned. The algorithm was tested using 12 instances of the pharmaceutical R\&D pipeline management problem. Solution qualities for the algorithm were 10\% better than the dynamic programming equivalent for all instances. The authors concluded that the algorithm was, nevertheless, computationally expensive when applied to the pharmaceutical R\&D pipeline management problems. \cite{RN20} proposed a knapsack-problem based heuristic (KDA) for pharmaceutical R\&D pipeline management problems under endogenous uncertainties. The KDA generates and solves a series of knapsack problems based on realized outcomes and has been shown to obtain high-quality feasible solutions for problems with up to 1,048,576 scenarios. \cite{RN11} developed a relaxed knapsack-problem based decomposition heuristic (RKDA) for the same problem to generate a valid dual bound. The relative gaps between KDA (primal bound) and RKDA (dual bound) were within 4.66\% for problems with over 1 million scenarios. \cite{RN26} further extended the KDA with different sub-problem generation approaches: every time period generation scheme (ETP) and at each realization generation scheme (AER), and further tested the KDA for various instances of the problem. \cite{RN27} developed a generalized knapsack-problem based decomposition algorithm (GKDA) to obtain feasible solutions for MSSPs under both endogenous and exogenous uncertainties. The GKDA is a scenario-free decomposition framework, which incorporates the initial and conditional NACs implicitly and addresses both space and time complexities of solving MSSP problems. The GDKA was applied to four planning problems and yielded feasible solutions with 0\% to 20\% relative gap for all cases.

In general, all approaches developed to date would realize significant computational benefits from reductions in numbers of NACs and scenarios as evidenced by heuristic approaches that employ scenario-free structures. 
 A full set of NACs can be generated by relating each scenario to every other scenario. This produces a very large number of constraints which are often redundant. Several authors have presented properties of MSSPs, which reduce the number of redundant NACs. Applicability of a lot of these properties is limited to MSSPs where the scenario set is considered complete (i.e., given by the Cartesian product) \citep{Goel-2006}. 

Often situations arise where the set of scenarios considered in the MSSP is given by a sub-set of the full set of scenarios. For instance, every scenario in the full set generated by the Cartesian product may not be realizable. Alternatively, a sub-set of the full scenario set may be used in cases where solving the MSSP with the full scenario set is computationally prohibitive. Properties used to generate NACs that are applicable for the full set of scenarios do not necessarily apply to sub-sets of scenarios. Thus, recent advances have focused on the generation of NAC sets for sets of scenarios which are not generated by the Cartesian product. \citet{Hooshmand-2016} developed two approaches; the first one used a mixed integer linear programming formulation (MILP) to find the minimum cardinality NAC set, and the second one employed graph theory, which added and removed NACs until there were no violations of non-anticipativity.  \citet{apap-2017} extended the graph-theory-based approach to handle both endogenous and exogenous uncertainty. \citet{boland-2016} developed a series of proofs which define the conditions required to omit NACs from the problem formulation. In all these works, the authors used the concept of a differentiator set, i.e., a set of sets, which store information on the realization of uncertainty required for two scenarios to be  distinguishable. The differentiator set is sufficient if the realization of the uncertain parameter is instantaneous, implying that the underlying value is revealed using a single indicator. If the realization of uncertainty is gradual there may be several indicators which result in partial revelation of uncertainty. The approaches developed to date cannot be applied to generate NACs of incomplete scenario sets for this type of uncertainty unless the underlying model of uncertain parameters is altered first.

\section{Problem Formulations and Background}\label{sec:form}

\subsection{MSSPs with endogenous and exogenous uncertainty}\label{sec:enex}
The presentation of the general form of the MSSP is adapted from \citet{apap-2017}, \citet{Hooshmand-2016} and \citep{Goel-2006}. Suppose that we have an MSSP where the planning horizon is defined as $\mathcal{T}=\{1,2,3,..,T\}$. 
 We first define a set of sources of endogenous uncertainties, $i\in \Ibb$, and a vector of uncertain parameters, $\theta_i$, associated with uncertainty source $i$. The realizable values of the uncertain parameter, $\theta_i$, are represented using a finite set $\Theta_i=\{\theta_i^1,\theta_i^2,\ldots,\theta_i^{\mathcal{R}_i}\}$. 
Further, we define $\Jbb$ as the set of exogenous uncertain parameters. We denote by $\xi_t$  a vector of exogenous realizations defined for parameters $j\in \Jbb$ and realized at time $t$. Each exogenous uncertain parameter is assumed to have a finite set of possible realizations, $\Xi_{jt}$.

The full scenario set is constructed using the Cartesian product, namely, the scenario set is defined as, $\mathbb{S} := \prod_{i \in I} \Theta_i \times \prod_{(j,t) \in \Jbb \times \T } \Xi_{jt}$.  
%The number of scenarios can be calculated using the cardinality of $\Theta_i$, where $|\mathbb{S}| = \Pi_{i \in I} |\Theta_i|$. 
We concentrate on the cases where the actual set of scenarios is a subset of the full Cartesian product, and hence we also define the corresponding set $S \subset \Sbb$.

\orange{
Further, we are primarily interested in the case when uncertain parameters allow for gradual realization. Intuitively, it means that as time progresses, depending on the decisions made, the set of possible values for each uncertain parameter shrinks, until the true values are eventually realized. 
This then means that for each endogenous parameter $i$ there exists a subset of $2^{\Theta_i}$ which defines the permissible states of realized information about it. 
We will refer to transitions between these states as \emph{events}. If all subsets of $2^{\Theta_i}$ are permissible, for each realization $\theta_i^k$ we can define an event $e_i^k$, indicating that $\theta_i^k$ has either been observed as the true realization of $\theta_i$ or eliminated.  In general, events correspond to the transitions between information states and can be decoded in different ways, for example, see the manufacturing example described below, where events correspond to completion of processing stages.
We can then define a set $\Ccal$ as the set of permissible combinations of events, corresponding to possible trajectories of uncertainty realization.  In the case of instantaneous uncertainty realization (or exogenous uncertainty) each parameter only has one event associated with it,  indicating that with a single event the underlying value of the uncertain parameter is realized.
}

% Further, we are primarily interested in the case when uncertain parameters allow for gradual realization. Intuitively, it means that as time progresses, depending on the decisions made, the set of possible values for each uncertain parameter shrinks, until the true values are eventually realized. Mathematically we will say that for each uncertain parameter, $\sigma  \in \Ibb \cup \Jbb$, there exist a set of events, $\mathcal{E}_\sigma$, which represents the events that must occur to fully realize $\sigma$. 
% In the case of instantaneous uncertainty realization (or exogenous uncertainty) $|\Ecal_\sigma |=1$ indicating that with a single event the underlying value of the uncertain parameter is realized. Further, depending on what kind of combinations of events are permissible, we can define a set of (ordered) subsets  $C \subset 2^\Ecal$ as the set of permissible events. For example, if the events must occur in a specific order, then not all subsets of events are realizable.

%Further, the gradual realization may be organized in a specific order, and hence the corresponding 
%if the gradual realization of uncertainty requires the set of events to be completed in a specific order similar to the manufacturing example, $\Ecal_\sigma$ is defined as an ordered set
%, \green{ $\Ecal_\sigma := \cap_{k\in [E]} e_\sigma^k$} .
%$\Ecal_\sigma = \{e_\sigma^1, \ldots, e_\sigma^{K_\sigma}  \}$.

At every time period the decision maker has to select values for decision  and the corresponding recourse variables. Observe that both should depend on scenario (i.e., the realization of uncertainty) and the set of realized events (i.e., the information available at the time of making the decision). We will then denote as $\bb_{tsc}$ and $\y_{tsc}$ the decision and recourse variables at time $t \in \T$ under scenario $s\in S$ given realized events $c \in C$. Note that if uncertainties are realized instantaneously, there is no need to  index on $c$. On the other hand, in the most general case of $C$ equal to the power-set over all realizations, the number of required variables can be prohibitively large. 

NACs must be enforced in the optimization formulation of MSSP, in order to ensure that the decisions do not take into account information that is revealed in the future. We  write the NACs in a general form as in \eqref{eq:genNAC}.
\begin{subequations}\label{eq:genNAC}
	\begin{align}
	& [Z_{t}^{rs}] \bigvee 
	\begin{bmatrix}
	\neg Z_{t}^{rs}\\
	\y_{tsc} = \y_{trc}\\
	\bb_{tsc} = \bb_{trc}
	\end{bmatrix}
	\quad \quad \forall r,s \in \Sbb, r\neq s, \forall t \in \T, \forall c \in C \label{eq:genNAC-disj}\\
	&Z_{t}^{rs} =  F(t, \bb, \y), \quad \quad \forall r,s \in \Sbb, \forall t\in\T \label{eq:genZ}
	\end{align}
\end{subequations}
Disjunction \eqref{eq:genNAC-disj} provides the NACs, given that $Z_{t}^{rs}$ is defined as the indicator variable for whether scenarios $s$ and $r$ are distinguishable at time $t$, and is defined in equation \eqref{eq:genZ}. Function $F$ is problem specific, and in \eqref{eq:genZ} $\bb$ and $\y$ refer to the whole collection of decision and recourse variables. %\red{I am a little confused on whether $Z$ should depend on $c$.}

Observe that NACs in \eqref{eq:genNAC} are written for all pairs of scenarios and all time periods. Depending on what the difference between  the scenarios  in a specific pair is, it may be possible to reduce the number of such constraints. Let us define $\tau(r,s)$ as the latest time period when scenarios $r,s$ are indistinguishable in exogenous uncertainty, i.e.,  $\tau(r,s)=\max⁡\{t\,|\, \xi_{t'}^s = \xi_{t'}^{r}  \,\, \forall t' \le t \}$. Then if two scenarios differ in exogenous uncertainty only, then \eqref{eq:genNAC-disj} reduces to 
$$	\begin{bmatrix}
%\neg Z_{tc}^{rs}\\
\y_{tsc} = \y_{trc}\\
\bb_{tsc} = \bb_{trc},
\end{bmatrix}, \quad \forall t \le \tau(r,s), \forall c \in C $$
and there is no need for disjunction, or equivalently, $Z_{t}^{rs} = 0$ for $t\le \tau(r,s)$, and $Z_{t}^{rs} = 1$ otherwise.
If two scenarios differ in both exogenous and endogenous uncertainty, then   \eqref{eq:genNAC-disj} can be written as 
$$[Z_{t}^{rs}] \bigvee 
\begin{bmatrix}
	\neg Z_{t}^{rs}\\
	\y_{tsc} = \y_{trc}\\
	\bb_{tsc} = \bb_{trc}
\end{bmatrix}
\quad \quad \forall r,s 
\in \Sbb, r\neq s, \forall t \le \tau(r,s), \forall c \in C .
$$

\orange{Note that this construction can also be viewed from the point of view of Markov decision processes (MDP) with discrete state-space. At any point in time, the pair $r$ and $c$, i.e., scenario and realized event set, determines the information available to the decision maker on the realization of uncertainty, and hence corresponds to the states (i.e., contains the information necessary to make decisions). NACs along with indexing of the variables with $c$ and $r$ corresponds to the idea that an optimal solution to an MDP is a \emph{policy}, i.e., a function prescribing actions to states. In the MDP framework NACs are usually not explicitly enforced, since the problem generally is not explicitly written in its entirety, and instead specialized algorithms are used.
}

In the next section we illustrate these constructions with a simple example, which we will then use to simplify the discussion of the proposed algorithm and as the basis for the numerical experiments.

%It is intuitively clear, that due to transitivity of equality, some of the NACs may be redundant. 

\subsection{An Illustrative Example}\label{sec:linear-representation}

 In order to illuminate the function of the endogenous NACs, specifically endogenous NACs with gradual realizations, here we present a simple demonstration, which  we will refer to as the manufacturing example. We will then use this example throughout the paper to  illustrate the developed algorithm. The example closely follows the description initially given in \cite{colvin-2008}.
 
% \red{Perhaps remove the detailed discussion of the example from here, and only introduce basics in the algo}

Suppose that a company produces two different products, $p\in\{P1,P2\}$. Both products are required to complete three ordered processing stages using the same equipment. A product may develop a defect during any of the processing steps. If a product develops a defect, it does not complete the next step. Whether or not a product completes all processing stages without a defect is not known with certainty. We can represent this uncertainty associated with product $p$ using an uncertain parameter $\theta_p$. The possible outcomes for the uncertain parameter $\theta_p$ are $\Theta_p=\{\omega_p^1,\omega_p^2,\omega_p^3,\omega_p^4\}$ where $\omega_p^1, \omega_p^2$, and $\omega_p^3$ represent a defect developed on product $p$ in processing stages one, two, and three, respectively. The outcome $\omega_p^4$  represents the case where the product $p$ was successfully manufactured with no defects.  In this example, the uncertainty is realized gradually, completing a processing stage will reveal if a defect is developed in that stage (a partial realization of uncertainty) but does not necessarily reveal whether or not the product will successfully complete all processing stages (complete realization of uncertainty). To have complete realization of uncertainty, all processing stages must be completed.

To illustrate how the non-anticipativity can be enforced for this example with the disjunctions \eqref{eq:genNAC-disj}, we begin by defining a binary variable, $\chi_{p,sg,t}^s$, which takes a value of 1 indicating that processing stage $sg\in[SG]$ of product $p\in P $ is started at time period $t\in \T $ for scenario $s\in \Sbb.$ 
Naturally it follows that we define a binary variable, $\zeta_{p,sg,t}^s$, which indicates whether the processing stage $sg$ has been completed for product $p$ at time $t$ in scenario $s$. For simplicity, in our example, we assume that each processing stage has a known duration for completion. We then can define intuitive
constraint $\zeta_{p,sg,t}^s=\zeta_{p,sg,t-1}^s + \chi_{p,sg,t-\delta_{p,sg}}$.
%, which relates the decision variable $\chi_{p,sg,t}^s$ to the indicator variable $\zeta_{p,sg,t}^s$ using the duration $\delta_{p,sg}$ required to complete the processing stage $sg$ for product $p$.
%\begin{align} \label{eq:zeta-chi}
%\zeta_{p,sg,t}^s=\zeta_{p,sg,t-1}^s + \chi_{p,sg,t-\delta_{p,sg}}
%\end{align}
Depending on the objective of the problem, the model may contain additional constraints or variables. These variables and constraints do not impact the development of NACs, as such, we omit them from this discussion.  Prior to completing any processing stages, all scenarios are indistinguishable. This implies that the values for the decision variable, $\chi_{p,sg,t}^s$, are identical until processing stages are completed.

For every scenario pair, $r,s \in \Sbb = \times_{p\in P} \{\theta_p \}$, there exists one or more differentiating event. In the case of gradually realized uncertainty, the differentiating event is an event which causes a partial realization of an uncertain parameter.  For example, consider two scenarios with the uncertain parameter realizations $(\omega_{P1}^1,\omega_{P2}^3 ) $ and $(\omega_{P1}^2,\omega_{P2}^4)$. The first scenario represents the case where product $P1$ develops a defect after processing stage one and product $P2$ develops a defect after processing stage three. In the second scenario, product $P1$ develops a defect after processing stage two and product $P2$ successfully completes all processing stages without developing a defect. Before the products complete any processing stages the scenarios are not  distinguishable. Once product $P1$ completes the first processing stage, these two scenarios become  distinguishable. If product $P1$ develops a defect during processing stage one, the underlying value for $\omega_{P1}$ is realized to be $\omega_{P1}^1$. If no defect is developed the value of $\omega_{P1}$ is realized not to be $\omega_{P1}^1$. To fully realize the value of $\omega_{P1}$, $P1$ must complete more processing stages.  

%\red{\emph{Word ``differentiable'', might be confusing due its calculus meaning. Are there any good substitutes, or is it a common way to say it in this area?}}

\orange{
In this example, the events that determine the gradual realization of uncertainty correspond to stage completion, hence we can define the set of events as $e^k_i$ for all $i,k$ as representing product $i$ completing stage $k$. Note that since the stages are completed in order the set of all permissible events can be defined as  $\Ccal = \cup_{p,q = 0}^3\{e^{P1}_1, \ldots, e^{P1}_q, e^{P2}_1, \ldots, e^{P2}_p\}$, where with slight abuse of notation we understand that $p = 0$ (or $q=0$) corresponds to no events occurring for product 1 (or 2). Note also that since the products are independent, the variables $\chi$ and $\zeta$ defined above only need to be indexed on the events for the corresponding product.
For any two scenarios $(r,s)$ and any set of events, we can establish whether the scenarios are distinguishable given the occurred events. We can then define $\Psi(r,s)$ representing the set of differentiating events that can occur causing scenarios $r$ and $s$ to become  distinguishable. For example, if scenario $r$ corresponds to both products passing all stages; and scenario $s$ represents $P1$ failing at stage 2, and $P2$ failing at stage $3$, then the set of differentiating events consists of P1 failing at stage 2, and P2 failing at stage 3, i.e., we can write $\Psi(r,s) = \{ e^2_{P1}, e^3_{P2} \}$.  We can then use $\Psi$ to define $Z_{rs}^t$ and the necessary disjunction based on variables $\zeta$, as in \eqref{eq:genNAC}. The disjunction can be converted to a linear constraint using standard integer programming techniques. 
}
\section{The Proposed Approach for Generation of NACs for Scenario Subsets}\label{sec:method}

\subsection{Graph-based model for NACs}\label{sec:graphModel}

Generation of NACs relies on knowledge of the scenario structure and how the groupings of indistinguishable scenarios change with the realization of uncertainty. 
%The goal of this work is to find an efficient approach to generate a minimum cardinality set of NACs, which we will denote by $\Ncal$. 
The algorithm proposed in this section will focus on minimizing the number of pairs of scenarios on which the NACs are enforced. Similarly, when discussing the reduction of the number of NACs we will compare the results in terms of the number of pairs of scenarios. Note that this is a standard approach in the literature, (see, for example \cite{Hooshmand-2016}). %Further, it is intuitively clear that reducing the number of pairs of scenarios for which NACs are applied also reduces the number of actual NACs.

 In order to construct the algorithm, we will consider a graph $\mathcal{G}=(\mathcal{S},E)$ with vertices representing scenarios $s\in \mathcal{S} \subset \Sbb$ and the edges   representing   NACs which relate the pairs of scenarios. \orange{Hence, a complete graph  would correspond to the basic formulation, where NACs are enforced between any two scenarios. A graph with just some arcs, corresponds to a formulation with only some NACs enforced. }
%	.  Here we assume that vertices represent scenarios $s\in \mathcal{S} \subset \Sbb$ and the edges will  represent  conditional NACs which relate the pairs of scenarios.  
	The graph can be naturally visualized by placing the vertices (scenarios) in a $|\Ibb| + |\Jbb|$ dimensional grid, according to the corresponding realizations of the uncertain parameters.  Before discussing the algorithm, we will establish some properties of the graph and how the grouping of vertices on the graph changes with the changes in the knowledge of the underlying values of the uncertain parameters in the system.

\begin{figure}
%	\FIGURE
\centering
	\includegraphics[width=0.6\textwidth]{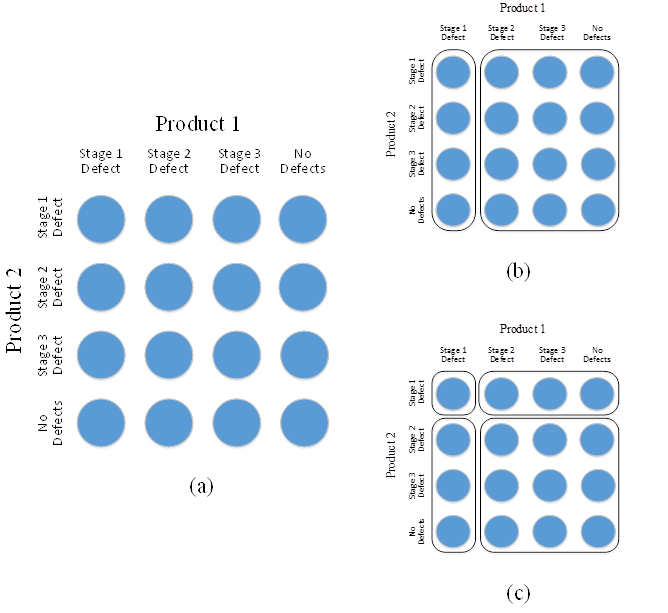}
	\caption{(A) The graph corresponding to the full set of scenarios in the manufacturing example arranged in a grid, without any NACs enforced. (B) Scenario groups of the manufacturing example formed by event set $\{e_{P1}^1\}$. (C) Scenario groups formed by the event set $\{e_{P1}^1,e_{P2}^1\} $\label{fig:manuf-full-set}} 
\end{figure}

Consider the set of permissible event subsets $C$, defined earlier. Any event set $c \in C$ defines a partition of the scenario set $S$, 
\begin{align*}
\Pi(c) : S = \bigcup_j S_j^c,
\end{align*}
where two scenarios $r, s$ are in the same set $S_j^c$ if the events in $c$ do not differentiate between $r$ and $s$. In the remainder of this paper we will refer to such subsets as \emph{cuts}.
Figure \ref{fig:manuf-full-set} illustrates this construction for the manufacturing example discussed earlier. Specifically, Figure \ref{fig:manuf-full-set}a shows the full scenario set arranged in a  grid according to realizations of uncertainty, while Figures \ref{fig:manuf-full-set}b and \ref{fig:manuf-full-set}c give the  partitionings for two event sets.

%Illustrating this concept with the manufacturing example, we begin by constructing the graph placed on a two-dimensional grid. Here we show the complete scenario set, however the concept can be generalized to any subset of the full scenario set. Figure \ref{fig:manuf-full-set}a shows the scenarios projected onto the lattice structure for the manufacturing example. Recall that for each product there are four possible outcomes, resulting in a total of 16 scenarios.

\begin{figure}
%	\FIGURE
\centering
	\includegraphics[height=0.9\textheight]{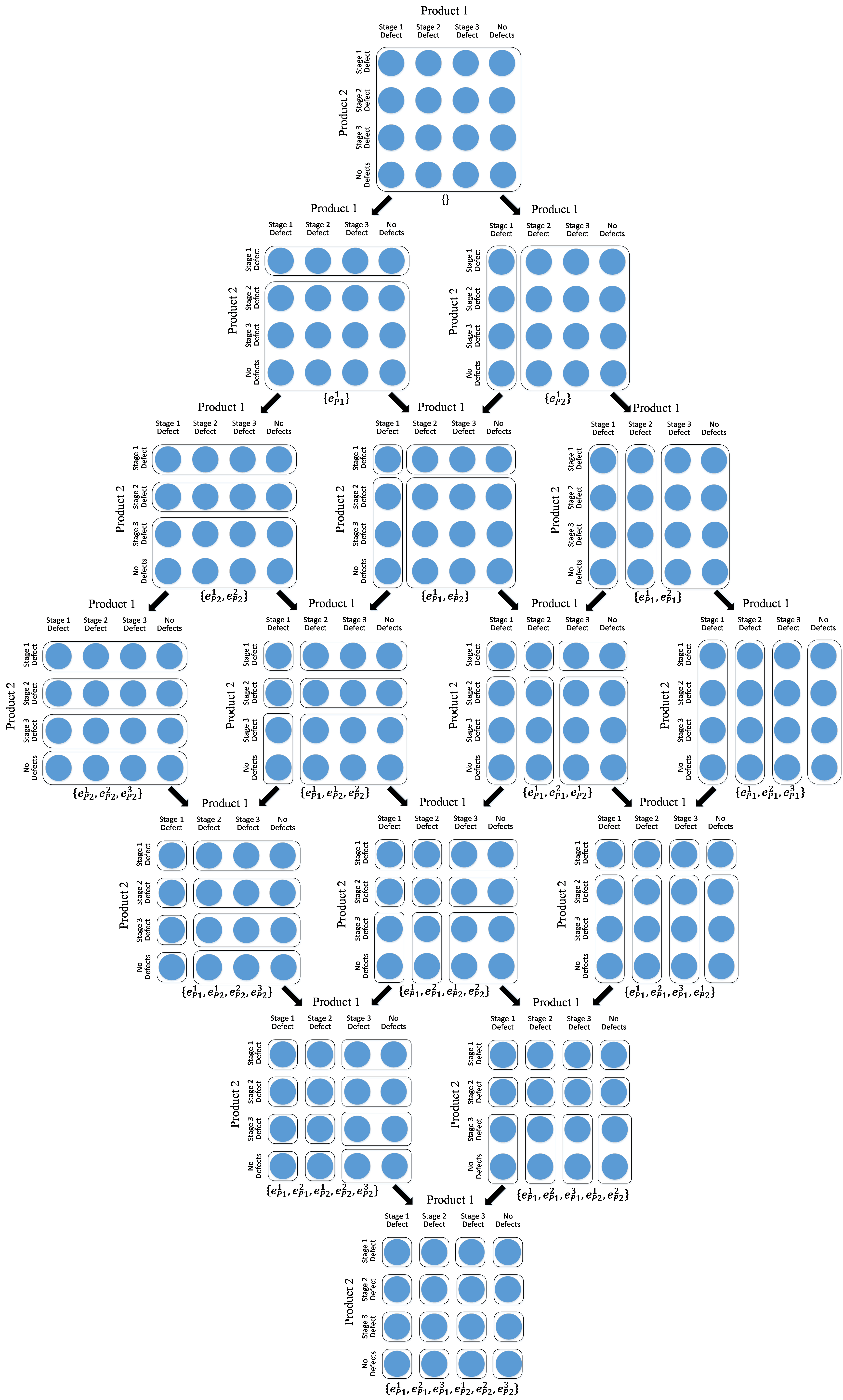}
	\caption{The sets formed by each event set $c\in \Ccal $ for the manufacturing example. \label{fig:ex-all-scen}} 
\end{figure}

Next, we will describe the relationship between a set of constraints sufficient to enforce nonanticipativity and connectivity of subsets of nodes. At the initial time period no uncertainty has been realized. This implies that all scenarios are identical, thus there is one group of scenarios (the whole set). Throughout the planning horizon uncertainty is realized. As uncertainty is realized, the scenario set is divided into subsets based on the realizations. Figure \ref{fig:ex-all-scen} shows the progression of uncertainty realization in the manufacturing example. Throughout the paper we will refer to the cardinality of event set $c$ as the \emph{degree of uncertainty realization}, since  $|c|$ represents the number of differentiating events which have occurred in order to reach the state of uncertainty realization. 
%Notice that the top set has $|c|=0$, and that we say that the degree of uncertainty realization $k=0$ for the initial set. Similarly, the set shown at the bottom of Figure \ref{fig:ex-all-scen}, where each scenario is in its own group, represents a degree of realization equal to six.

%\orange{\begin{definition}
%	We will say that \emph{the degree of uncertainty realization} of a cut set $c\in \Ccal$ is given by $|c|$, or in other words,  $|c|$ represents the number of differentiating events which have occurred in order to reach the state of uncertainty realization.
%\end{definition}}

%Increasing the degree of uncertainty realization corresponds to the occurrence of a differentiating event.
\begin{lemma}\label{prop:part-lemma} % \red{rephrase according to reviewer}
If $c\subset c'$, then partition $\Pi(c')$ is a refinement of $\Pi(c)$, that is for every $j'$ there exists $j$ such that $S_{j'}^{c'} \subset S_j^c$.
\end{lemma}
\begin{prf}
	The claim is a direct consequence of construction of partitioning. Indeed, $c \subset c'$ implies that the state of the system corresponding to events in $c'$ can be achieved as an evolution of the system from a state corresponding to events in $c$. Consequently, scenarios that are indistinguishable in this future state had to be indistinguishable in the past as well. 
%	\red{update notation}
%The order of the cut set refers to the number of differentiating events which have occurred to reach the state of uncertainty realization. 
%%The scenarios $s\in \Sbb $ form subsets $\sfr\in \Sf_c$ based on differentiating events. More specifically, 
%Two scenarios $r,s\in \sfr\in \Sf_c$ if they are indistinguishable with respect to the differentiating events in $c$. Increasing the number of differentiating events in the cut set $c$ (i.e. $|c' |=|c|+1$ and $c\subset c'$) results in formation of a new set of subsets $\Sf_{c' } $. 
%If two scenarios $\{s,s' \} \notin \sfr\in \Sf_c $ then the event differenting scenarios $s$ and $s'$ exists in the set $c$. This implies that if an additional differentiating event is added to the cut set $c$ to form a cut set $c'$, then $\{s,s' \}\notin \sfr\in \Sf_{c' }$.
%\Halmos 
\end{prf}

If two scenarios $r,s\in S^c_j$ for some cut $c\in \Ccal$, then they are indistinguishable for some permissible state of the system, and hence, non-anticipativity should be enforced for these two scenarios. On the other hand, a direct non-anticipativity constraint for $(r,s)$ can be eliminated if another path $\mathcal{P}$ exists within nodes in $S_j^c$ connecting $r$ and $s$ due to transitivity of the constraints. Further, since while uncertainty is gradually realized new scenario sets are constructed as subsets of existing sets (Lemma \ref{prop:part-lemma}), the sets formed by the differentiating event must exist as connected sets prior to the occurrence of the differentiating event. This observation leads to the following characterization of relationship between NACs and the corresponding graph.  

% To enforce the proper NACs, the resulting scenario sets created by the occurrence of the differentiating event must form a connected graph (Property 1). 

\begin{lemma}\label{prop:nac-sufficency}
	A set of constraints is sufficient to enforce nonanticipativity, if and only if for all cuts $c\in \Ccal$ all of the corresponding scenario subsets $S^c_j$ are connected in the corresponding induced subgraph.
\end{lemma}
\begin{prf}
Follows by construction, from the discussion above.
%\Halmos
 \end{prf}
%Proposition \ref{prop:nac-sufficency} establishes the structure of graph $\mathcal{G}$ that is required in order to enforce NACs. 
Next, we will demonstrate how
%
%
%This implies that the addition of at most $m-1$ edges is required to combine $m$ subsets when decreasing the order of uncertainty realization from $|c|$ to $|c|-1$.
%
using the knowledge of the structure of the subsets of vertices formed by each event set and the relationship between event sets, it is possible to develop an algorithm which generates a set of edges for the graph G such that it satisfies Lemma \ref{prop:nac-sufficency} with the minimal number of corresponding pairs of scenarios.
%\red{\emph{Minimal number of edges or number of NACs? Is there a difference}}

%\red{\emph{I don't think second part is needed}}
%\green{\begin{definition}
%	Any edge $(r,s)\in E$ is defined as \emph{necessary} if two conditions are true:
%	\begin{enumerate}
%		\item 	$\{r,s\}\in \sfr\in \Sf_c$ and for all $c'\in C$ if $|c|<|c' |$ then  $\{r,s\} \notin s'$ for all $s'\in \Sf_{c' }$  and,
%		\item second property
%	\end{enumerate}
%\end{definition}
%}
%	If there exists no edges (r^',s^' )∈E for which condition (1) holds 

\begin{definition}\label{def:necess}
	An edge $(r,s) \in E$ will be called \emph{necessary} for a set of nodes $S^c_j \in \Pi(c)$ if $r,s \in S^c_j$,  and $r,s \notin S^{c'}_i$ for all $S^{c'}_i\in \Pi(c')$ and $c'\in C$ such that 
	$|c|<|c'|$.
\end{definition}
%
%
%\green{	
%	The first condition in Definition \ref{def:necess} establishes that an edge $(r,s)$ is necessary to ensure Property 1 applies to the subsets formed by the cut set $c$. The premise is that if two scenarios $r$ and $s$ exist in a subset $s$ in the set of subsets $\Sf_c$ then the edge $(r,s)$ is necessary to ensure Property 1 if there exists no edge $(r',s' )\in E $ which already ensures Property 1 holds. It is also required that the edge $(r,s)$ in question must connect two scenarios which have not been present in any subset $s'\in \Sf_(c' )$ for all $c'\in C$ such that $|c' |>|c|$. It is important to note that there is no condition that requires $|c' \cap c|=1$. The lack of this condition implies that any edge defined as necessary in subset created by a cut set which has a higher order can be used to satisfy Property 1 in lieu of adding a new edge.
%}

 By this definition, an edge $(r,s)$ is necessary for $S^c_j \in \Pi(c)$ if non-anticipativity for scenarios $r,s$ is required (as $r,s\in S^c_j$), yet it cannot be established based on edges that are within subsets in $\Pi(c')$ for $|c'| > |c|$, i.e., when more information on uncertainty realization is available.   In other words, unless non-anticipativity for $r,s$ is enforced with a path within $S_j^c$, then it will not follow from any other set of NACs. 
% Note that it is not required in the definition above to have $|c'\setminus c| = 1$, as non-anticipativity can be enforced based on paths in $\mathcal{G}$ created based on a cut set which has a higher order  in lieu of adding a new edge.

This then leads to an intuitive idea for generating a set of NACs. The proposed algorithm is presented in Algorithm \ref{alg:alg}. It is constructed as a greedy procedure, which begins by considering the fully realized uncertainty set ($k = \max\{|c |, c\in \Ccal \}$) and then generates subsets of necessary edges separately for each possible event set $c \in C$ such that $|c| = k$. 
%After generating subsets, MSTs are determined for each subset assuming that the graph already contains the edges deemed necessary at the previous iterations. 
Any edge deemed necessary is stored, then the value of $k$ is decreased. The process is repeated until $k=0$. The resulting NAC set is represented by the edges deemed necessary in all iterations. 
\begin{algorithm} % enter the algorithm environment
	%	\caption{Calculate $y = x^n$} % give the algorithm a caption
	% and a label for \ref{} commands later in the document
	\begin{algorithmic}[1] % enter the algorithmic environment
		\STATE $k := \max\{|c |, c\in \Ccal \}$
		\STATE $\mathcal{N} := \emptyset$
		\WHILE{$k \ge 0$}
		\FOR{$c\in \Ccal : |c|=k$}
		\STATE \hspace{1cm} Generate subsets $S_j^c\in \Pi(c)$
		\STATE \hspace{1cm} Construct a spanning tree over the restriction of the graph onto nodes in $S_j^c$ and existing edges $\Ncal$ \label{alg:step}
		\STATE \hspace{1cm}  Add edges identified in STEP \ref{alg:step} to $\Ncal$		
		\ENDFOR
		\STATE $k=k-1$
		\ENDWHILE
	\end{algorithmic}
	\caption{The Sample Non-Anticipativity Constraint (SNAC) algorithm for minimum cardinality NAC generation \label{alg:alg}}
\end{algorithm}

An example progression of the algorithm is presented in Section \ref{sec:comp-simple}. 
\begin{theorem}\label{prop:alg-correctness}
	Algorithm \ref{alg:alg} terminates with a subset of edges which corresponds to a set of NACs enforced on the minimum number of pairs of scenarios. The algorithm can be implemented in a way that guarantees that the running time   grows as $O(|\Ccal| |S|^3)$ as $|\Ccal|, |S| \rightarrow +\infty$, where $S$ is the number of scenarios and $\Ccal$ is the the set of permissible event sets.
\end{theorem}
\begin{prf}

Correctness and optimality of the resulting NAC-set follow by construction and Lemmas \ref{prop:part-lemma} and \ref{prop:nac-sufficency}. Indeed, the algorithm explicitly enforces connectivity on subsets $S_j^c$ for all $c$ and $j$. 
\orange{Further, all arcs added at step \ref{alg:step} of the algorithm are necessary in the sense of definition \ref{def:necess} by construction. 
Hence it is impossible to remove any of the arcs without violating connectivity. Indeed, suppose we remove an arc connecting scenarios $r$ and $s$. In this case, the NACs for these two scenarios are not enforced since the arc was added as a part of the spanning tree, and these NACs could not have been enforced with any of the previously added arcs by construction. 
} 
Finally, while each component $S_j^c$  potentially can be   connected in multiple ways, this choice has no effect on the number of arcs needed to be included at future iterations. This claim follows from Lemma \ref{prop:part-lemma} and the way each subset $S_j^c$ is constructed by combining a number of smaller subsets $S_i^{c_i}$ with $|c_i|>|c|$, all of which are already connected since they have been considered at a previous iteration.

In order to establish running time estimate, observe that the algorithm iterates over all $c \in C$ and then over all $S_j^c \in \Pi(c)$. For each $ S_j^c$ it then generates the set of scenarios in $S_j^c$ and introduces arcs to enforce connectivity. The number of subsets $S_j^c$ for each $c$ is bounded by $|S|$. In a simplest implementation connectivity of a graph can be enforced with a greedy algorithm constructing a spanning tree in $O(n^2)$, where $n$ is the number of nodes. In our case, the number of nodes in each subset $S_j^c$ is bounded by $|S|$. Further, generation of each subset $S_j^c$ can be organized in $O(|S|)$ by enumeration. This then results in overall $O(|C||S|^3)$ running time.
	%consists of three main stages: outer loop iterating over all $c \in C$,  an inner loop iterating over all $\sfr \in \Sf_c$ followed, finally,  generation  of subset $\sfr$ and MST procedure.
\end{prf}

\begin{remark}
	
	Each component $S_j^c\in \Pi(c)$ for $|c| = k$ is treated ``in parallel'', i.e., spanning trees to enforce connectivity are constructed without regard to  edges that are added at the same iteration. This can lead to creation of cycles in $\mathcal{G}$. When connectivity is enforced on any subset, each connected component is treated as a single node.
\end{remark}

The running time is proportional to $O(|\Ccal|)$, which, in the worst case, is a powerset of the set of uncertainty realization events, i.e., it is exponential in the number of events. On the other hand, if the number of uncertain parameters and events is fixed, and $S \ll \mathbb{S}$, which is often the case, then the running time is polynomial in the number of scenarios.

\begin{remark}
Note that the description of the algorithm is provided in terms of event sets and does not emphasize whether 	underlying uncertainty is endogenous or exogenous. Indeed, while in the case of exogenous uncertainty, one can trivially avoid having to model gradual realizations, its can still be processed by SNAC in the exact same way. For an exogenous uncertain parameter (eg., demand in the manufacturing example) the corresponding  event set $\Ecal_\sigma$ consists of a single member, yet the algorithm can still be allied in the same manner as for purely endogenous cases. \orange{As an illustrative example, consider a case similar to the manufacturing example with single product, and instead of the uncertainty due to processing stages of product two, we consider uncertain demand. Suppose that demand parameter can take any of the values $\{D1, D2, D3, D4\}$, which are gradually realized in two phases: at time $t=1$ the demand is revealed to be either $\{D1, D2\}$ or $\{D3, D4\}$ and then at time $t=2$ the actual value is observed. Then, following the same construction as in the example above, we can denote events as $\{t=1, t=2,e_{p1}^1,e_{p1}^2,e_{p1}^3\} $ and then construct set $\Ccal$ and the diagram analogous to Figure \ref{fig:ex-all-scen} as given on Figure \ref{fig:exo-scen}. Consequently, the algorithm can be applied in the same way as for the purely endogenous case. 
}

\begin{figure}
%	\FIGURE
	\includegraphics[width=0.9\textwidth]{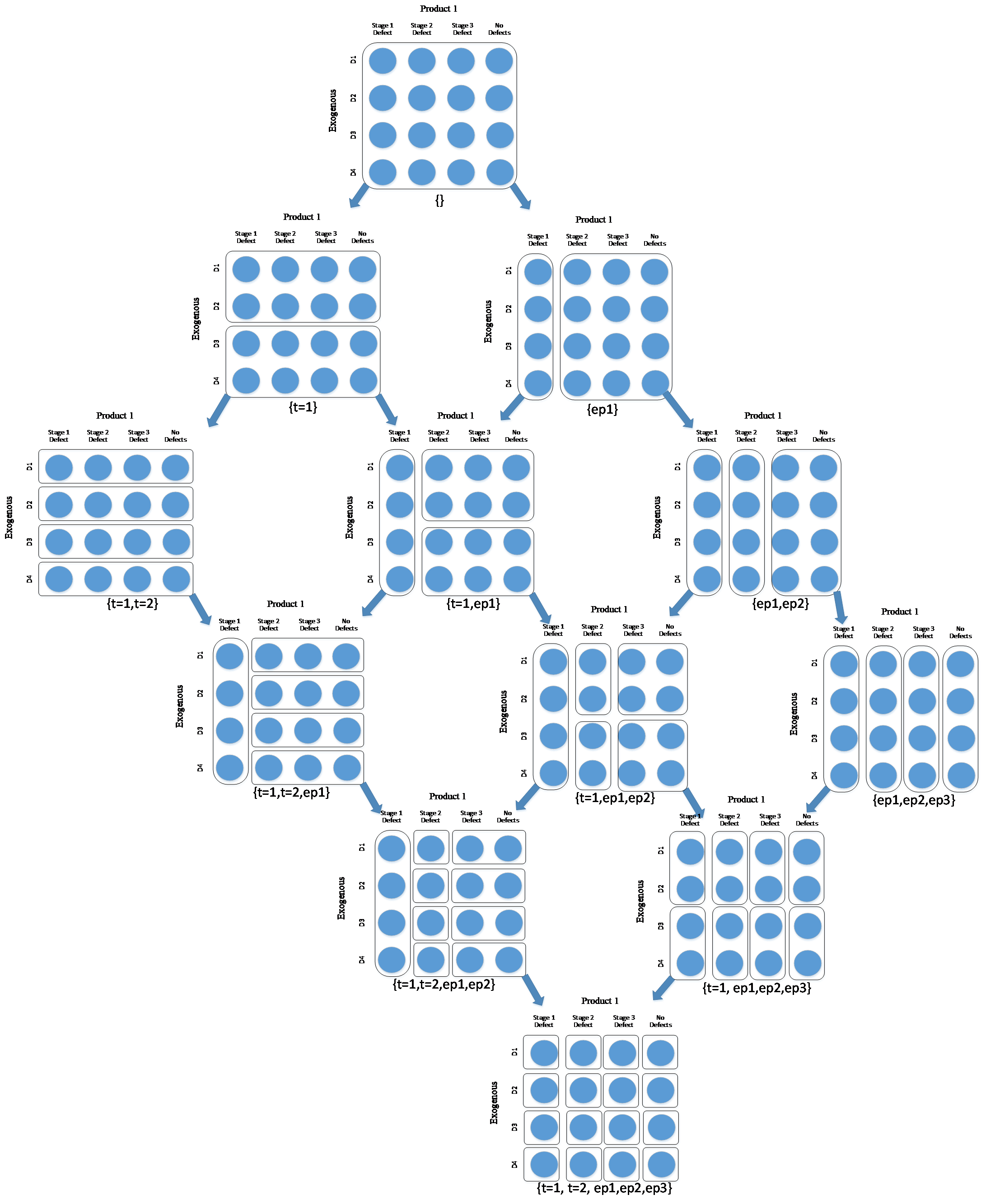}
	\caption{The sets formed by each event set $c\in \Ccal $ for the manufacturing example with exogenous demand uncertainty. \label{fig:exo-scen}} 
\end{figure}
\end{remark}

\subsection{Illustration of the Algorithm Progression for a Simple Manufacturing Example}\label{sec:comp-simple}

%\green{The manufacturing problem is characterized by a set of products which are required to complete a series of ordered processing steps. In our example, we assume that there are two products $p\in \{P1,P2\}$ which are required to complete three ordered processing steps $sg\in \{1,2,3\}$. In each processing step, the product may or may not develop a defect. We represent the uncertainty in whether each product develops a defect as an endogenous uncertain parameter. The parameter may take one of four discrete realizations, $\omega_p\in \{\omega_p^1,\omega_p^2,\omega_p^3,\omega_p^4\}$ where the numbers $\omega_p^{1-3}$ represent the failure at stages 1, 2 and 3. And $\omega_p^4$ represents a successful completion of all processing stages. The MSSP formulation will contain an objective function, likely an economic measure as well as scenario specific constraints. For example, constraints exist in the formulation which ensure that tasks within each scenario are completed in the correct order. Here, we limit our discussion to the set of NACs which ensure that the decisions do not anticipate the outcome of the uncertain parameters. }

\begin{figure}[H] \centering
%	\FIGURE
	\includegraphics[height=0.3\textwidth]{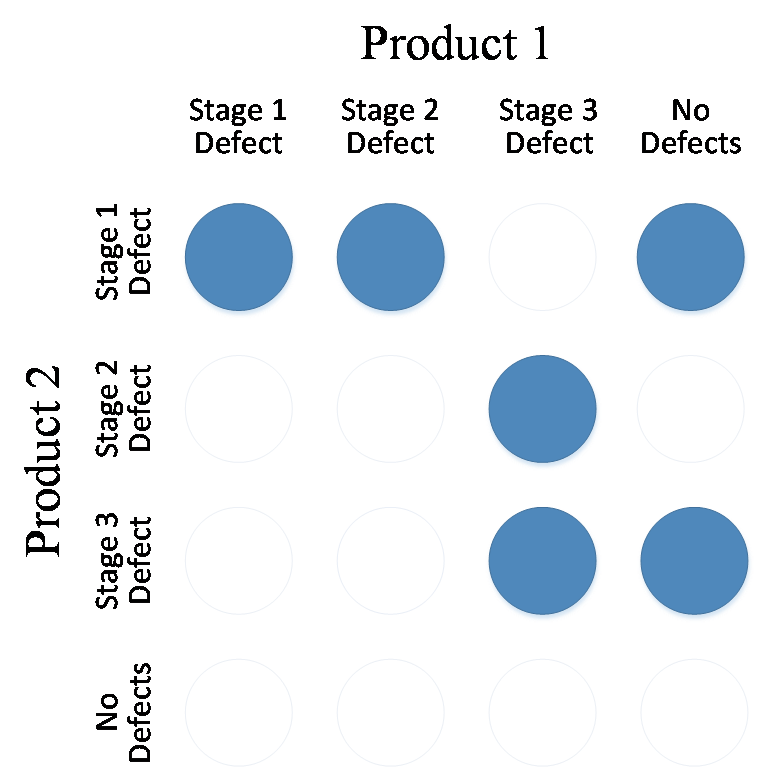}
	\caption{The six sampled scenarios. \label{fig:six-scen}} 
	\end{figure}

The manufacturing problem considered is described in detail in Section \ref{sec:linear-representation}. The cardinality of the full scenario set for the manufacturing problem with two products and three processing stages is 16 (calculated as $4^2=16$). For this example, we sample 6 random scenarios, $S=\{(\omega_{P1}^1,\omega_{P2}^1 ),(\omega_{P1}^4,\omega_{P2}^3 ),(\omega_{P1}^2,\omega_{P2}^1 ), (\omega_{P1}^3,\omega_{P2}^2 ), (\omega_{P1}^4,\omega_{P2}^1 ), (\omega_{P1}^3,\omega_{P2}^3 )\}$.
Figure \ref{fig:six-scen} gives the initial graph corresponding to these scenarios arranged in a grid.
 Table \ref{t:cuts} shows the event set $C$ (organized by $|c|$, referred to as \emph{order}).
\begin{table}[h!]
	\small
	\centering
	\caption{The event sets $\Ccal$ for the manufacturing example.}
	\label{t:cuts}
	\begin{tabular}{cl}
		\hline
		Order&	$\Ccal$	\\
		\hline
		0 & 	\{\}\\
		1	&$\{e_{P1}^1 \}, \{e_{P2}^1 \}$\\
		2	&$\{e_{P2}^1,e_{P2}^2 \},\{e_{P1}^1,e_{P2}^1 \},\{e_{P1}^1,e_{P1}^2 \}$\\
		3	&$\{e_{P2}^1,e_{P2}^2,e_{P2}^3 \}, \{e_{P1}^1,e_{P2}^1,e_{P2}^2 \}, \{e_{P1}^1,e_{P1}^2,e_{P2}^1 \},   \{e_{P1}^1,e_{P1}^2,e_{P1}^3 \}$\\
		4	&$\{e_{P1}^1,e_{P1}^2,e_{P2}^1,e_{P2}^2 \},\{e_{P1}^1,e_{P1}^2,e_{P1}^3,e_{P2}^1 \}
		\{e_{P1}^1,e_{P2}^1,e_{P2}^2,e_{P2}^3 \}$ \\
		5	&$\{e_{P1}^1,e_{P1}^2,e_{P2}^1,e_{P2}^2,e_{P2}^3 \},\{e_{P1}^1,e_{P1}^2,e_{P1}^3,e_{P2}^1,e_{P2}^2 \}$\\
		6	&$\{e_{P1}^1,e_{P1}^2,e_{P1}^3,e_{P2}^1,e_{P2}^2,e_{P2}^3 \}$\\
		
		\hline

	\end{tabular}
\end{table}

\begin{figure}[H] \centering
%	\FIGURE
	\includegraphics[height=0.3\textwidth]{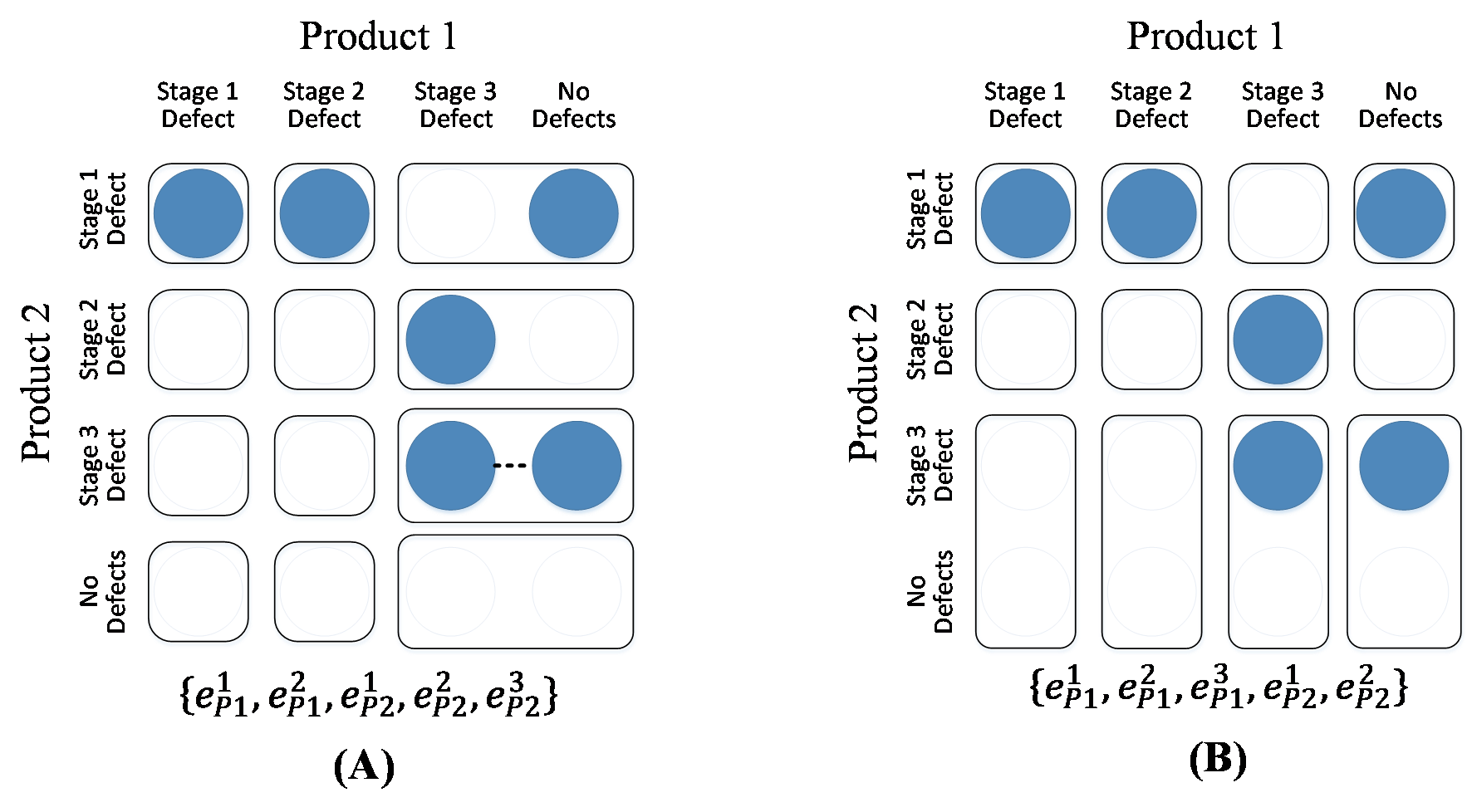}
	\caption{Subsets formed by the event sets of length five  \label{fig:alg-k-5}}
	
\end{figure}

\begin{figure}[H] \centering
%	\FIGURE
\centering
	\includegraphics[height=0.3\textwidth]{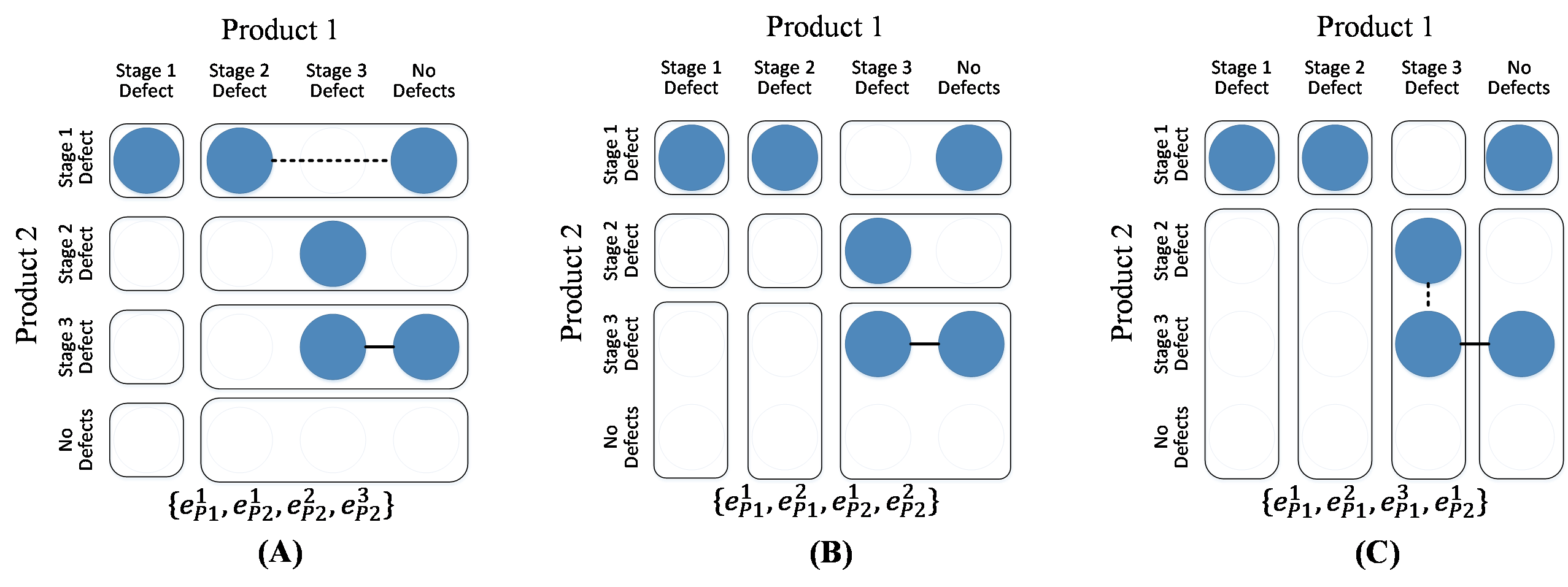}
	\caption{Subsets for the iteration where $k=4$  \label{fig:alg-k-4}}
	
\end{figure}

\begin{figure}
%	\FIGURE
	\includegraphics[width=\textwidth]{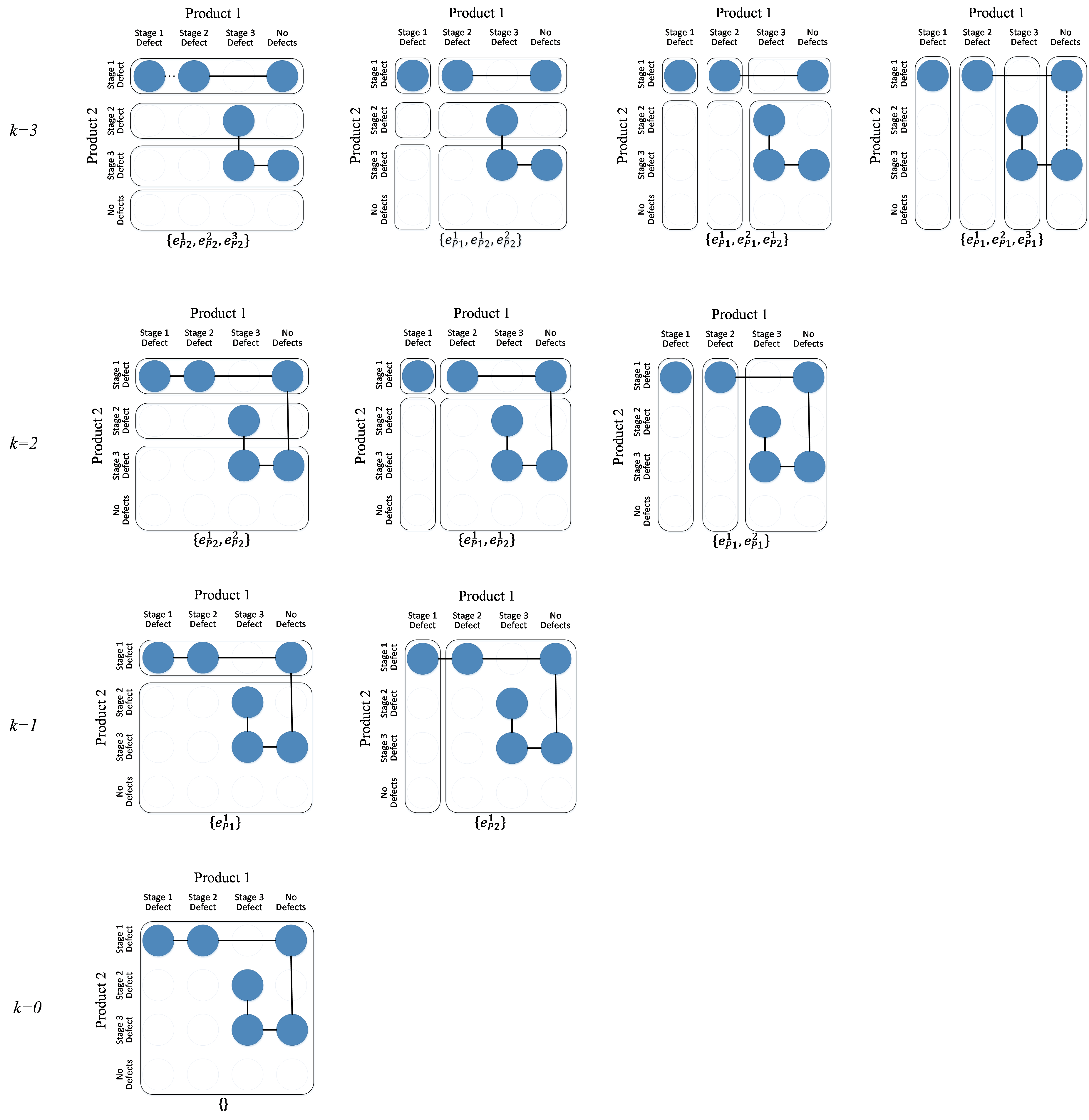}
	\caption{The last four iterations of the SNAC algorithm for the two product manufacturing example \label{fig:alg-end}} 
	
\end{figure}

The algorithm starts by initializing $k$ to six. The event set $c$ with $|c|=6$ produces subsets where each scenario forms its own subset. As a result, no edges are deemed necessary. Reducing $k$ by one and generating new subsets yields the groups given in Figure \ref{fig:alg-k-5}. The corresponding event set with five events is shown below each set of subsets. In Figure \ref{fig:alg-k-5}(A), two scenarios fall into the same subset
%. There exists no edge in the graph which can ensure that the corresponding components are connected. 
and the corresponding necessary edges  are represented with dotted lines. Figures \ref{fig:alg-k-4} and \ref{fig:alg-end} illustrate the remaining progression of the algorithm. It generated NACs for a total of five pairs  to connect six scenarios. Without reductions the number of pairs of scenarios with enforced NACs would have been 15 calculated by $0.5[|\Scal|(|\Scal|-1)]$.

%the value of $|\cup_{\sigma \in \{\omega_{P1},ω_{P2} \}} \Ecal_\sigma |$. In this case, the value of $k$ is set to six. The next step is to divide the scenarios into subsets for each $c\in \Ccal:|c|=k$. The cut set $c$ with $|c|=6$ produces subsets where each scenario forms its own subset. As a result, no edges are deemed necessary. Reducing $k$ by one and generating new subsets yields the groups given in Figure \ref{fig:alg-k-5}. The corresponding cut set with five events is shown below each set of subsets. In Figure \ref{fig:alg-k-5}(A), two scenarios fall into the same subset
%. There exists no edge in the graph which can ensure that the corresponding components are connected. 
%and MST provides the list of necessary edges represented with dotted line.

%After identifying the necessary edges, they are added to the set $\Ncal$, and the value of $k$ is decreased by one. Subsets are then generated for cut sets with four events. These sets are shown in Figure \ref{fig:alg-k-4}. The edges which were deemed necessary in the previous iteration are shown as solid lines in Figure \ref{fig:alg-k-4}. In Figure \ref{fig:alg-k-4}(B), the edges added during the previous iterations are sufficient to satisfy property 1 for all subsets. In both Figure \ref{fig:alg-k-4}(A) and \ref{fig:alg-k-4}(C), there exists at least one subset where additional edges are needed. The additional edges deemed necessary by the MST for each subset are shown with dotted lines, and they are added to the set $\Ncal$.

%The value of $k$ is decreased again, and iterations continue until the value of $k$ is equal to zero. The subsets and necessary edges for the remaining iterations are summarized in Figure \ref{fig:alg-end}. The iteration with $k=3$ adds the final two edges to the necessary set of edges. 
%\green{At this iteration, it is not possible to determine that all necessary edges are included in set $\Ncal$. The algorithm proceeds with the remaining iterations ($k=2$, $k=1$, and $k=0$), and determines that all generated subsets are connected with the existing edges.} 
%The algorithm generated a total of five NACs to connect six scenarios. Without reductions the number of NACs would have been 15 calculated by $0.5[|\Scal|(|\Scal|-1)]$.
% For this example, the SNAC algorithm reduced the total number of NACs by more than 50\%.  

\section{Computational Studies}\label{sec:comp}

%In this section, we first present a step by step example of the generation of the NAC set for the simple manufacturing example introduced in Section \ref{sec:linear-representation} using SNAC algorithm (Section \ref{sec:comp-simple}). The example considers two different products, each of which should complete three ordered processing stages.  

%In this section, we present results of two computational studies. First, we investigate SNAC algorithm itself in terms of the running time and the reduction in the number of scenario pairs requiring NACs. In the second study we solve some simple MSSPs to evaluate whether application of SNAC algorithm results in solution time savings.
%
%. The study compares the number of pair of scenarios for NACs generated by the SNAC algorithm to the number of NACs that will be added to the MSSP without any reductions. 
%The computational studies cover from two to six uncertain parameters with between two and ten realizations with varying number of scenarios. A detailed description of each case study can be seen in \red{Appendix A}. 
In this section we present results of a computational study. The aim of the study was to evaluate whether the proposed algorithm is able to significantly reduce the number of scenario pairs requiring enforcement of NACs and further, whether it then translates into savings in terms of MSSP solution time.   
All computational studies were implemented in Python Version 3.5.0 and performed on a 64-bit machine running 
Windows 7 with Xeon CPU E3-1241 and 32 GB RAM.  The resulting optimization problems were solved in IBM ILOG CPLEX.

MSSP instances 	solved were generated in the same way as the models presented in \cite{colvin-2008}. \orange{The problems follow the same logic as the manufacturing example we used above with more products, stages and scenarios. For completeness, we present the problem formulation and parameter values in \ref{sec:app_form} and \ref{sec:app_params} respectively. For a more detailed description of the cases, see   \cite{colvin-2008}.}
Note that here we use the same terminology as in the manufacturing example to describe the observed results. We vary the number of products, the number of possible outcomes (processing stages) and the number of scenarios sampled from the full Cartesian product.
	In each case, random instances are generated 30 times. \orange{In each instance generated, we solved both the original MSSPs without identifying the minimum scenario pairs (referred to as ``Full'') and the MSSPs constructed using the proposed algorithm (referred to as ``SNAC''). Tables \ref{t:MSSPstats} and \ref{t:MSSPtime} report the information on the instance sizes and running time respectively.} It should be noted that the solutions of both models yielded the same objective function value for all instances and scenario subsets.
	
\orange{First, consider Table \ref{t:MSSPstats}, which reports 	the number of NACs and the corresponding number of scenario pairs due to SNAC (average and maximum over 30 random instances) and compares it with the full problem. We can observe that both values quickly grow as problem size increases, yet in all cases the algorithm significantly reduces the number of NACs and scenario pairs. This result is not suprizing, since it is clear that full formulation contains a significant number of redundant constraints. Next we investigate whether reduction in NACs leads to a reduction in solution time.
	}  
	It should be noted that the SNAC only utilizes the uncertain parameter and outcome information without considering the planning horizon or potential completion times. The number of NACs for certain problems including these case studies can be further reduced if the task completion times and planning horizon is considered as shown in the NAC reduction properties by \cite{colvin2010}.
	
\orange{Table \ref{t:MSSPtime} reports running time. Specifically, we report three separate times: SNAC algorithm time, model generation time (i.e., the time required to generate NACs) and CPLEX time 	(time used by the solver once the model is loaded). Here we can observe a trade-off, whereby on one hand, in all cases the models created by SNAC take significantly less time to both generate and to obtain a solution; yet on the other hand, SNAC itself may require a significant time to complete. For example, in case of 5 products, 4 outcomes and all 1024 scenarios, the algorithm used on average 5140.6 seconds to run, which then allowed to reduce model generation time from 3454.9 to 229.4, and CPLEX time from 4331.7 to only 43.5. This can be very beneficial for two reasons. First, the total time with SNAC ($5140.6 + 229.4 + 43.5 = 5413.5$) is 30\% less than the total time for the full problem ($3454.9 +4331.7 = 7786.6$). Secondly, note that SNAC only uses scenario structure and is independent from the rest of model, so if a decision maker is interested in solving a collection of problems, which differ in parameter values (e.g., costs, revenues, demand, etc), they only need to run SNAC once, and all subsequent problems can be solved with the same optimal NAC structure. In the case of full problem formulation, each new problem would need to be solved from from scratch. }

\orange{
To sum up the case study, we conclude that our proposed algorithm can significantly reduce teh number of NACs that need to be enforced in an MSSP with gradual realization of uncertainty. Further, this reduction can directly result in computational time savings. These savings are especially significant if one needs to solve a collection of problems with the same scenario structure.
}

\begin{table}[h!]
	
	\centering
	\caption{Summary of problem sizes for case studies solving sample MSSPs based on the manufacturing example. All values are averages or maximums over 30 instances. Average and maximum for numbers of pairs of scenarios and actual NACs  are reported for SNAC algorithm. Columns Full report the number of NACs and scenarios pairs in the full problem (does not change with each random instance).	\label{t:MSSPstats}}

	\begin{tabular}{lll|lll|lll}
		%		\toprule
		
		\hline
		&&& \multicolumn{3}{c|}{\# of NACs} & \multicolumn{3}{c}{\# of scenario pairs}\\
		Products	&	Outcomes	&	Scenarios	&	Average &	Max	&	Full 	&	Average	&	Max	&	Full 
		\\
		\hline
		2	&	4	&	12	&	793.6	&	840	&	3192	&	16	&	17	&	66	\\
		2	&	10	&	24	&	4593.6	&	4944	&	39792	&	31.6	&	34	&	276	\\
% 		\hline
		3	&	4	&	6	&	1351.2	&	1998	&	2988	&	6.7	&	10	&	15	\\
		3	&	5	&	24	&	10869.6	&	13008	&	72936	&	40.9	&	49	&	276	\\
% 		\hline
		4	&	3	&	12	&	1730.7	&	2128	&	5328	&	21	&	26	&	66	\\
		4	&	4	&	128	&	39196	&	42032	&	975872	&	322.4	&	346	&	8128	\\
		4	&	5	&	24	&	9370.7	&	11136	&	44256	&	58	&	69	&	276	\\
% 		\hline
		5	&	4	&	64	&	37350	&	42170	&	302720	&	246.9	&	279	&	2016	\\
		5	&	4	&	1024	&	581120	&	581120	&	78571520	&	3840	&	3840	&	523776	\\
 		\hline
		%		\midrule
		
		%		\bottomrule
		%		
		
	\end{tabular}
\end{table}

\begin{table}[h!]
	\small
	\centering
	\caption{Summary of solution times for case studies solving sample MSSPs based on the manufacturing example. All times are reported in seconds and are averages over 30 instances. SNAC Running Time refers to the running time of the algorithm, 
		Model Generation is the time required to generate the NACs (either all or based on minimum pairs), and CPLEX represents the time used by CPLEX to obtain the optimal solution (either full problem or the one produced by SNAC). \label{t:MSSPtime}}
	
	\begin{tabular}{lll|l|ll|ll}
		%		\toprule
		\hline
		&&& SNAC Running& \multicolumn{2}{c|}{Model Generation} & \multicolumn{2}{c}{CPLEX }\\
		Products	&	Outcomes	&	Scenarios	&	Time	&	SNAC	&	Full	&	SNAC 	&	Full
		\\
		\hline
		2	&	4	&	12	&	0.321	&	0.123	&	0.215	&	0.142	&	0.234	\\
		2	&	10	&	24	&	1.597	&	0.546	&	2.317	&	0.648	&	2.337	\\
		3	&	4	&	6	&	0.707	&	0.205	&	0.266	&	0.245	&	0.284	\\
		3	&	5	&	24	&	1.295	&	1.101	&	4.2	&	1.261	&	5.147	\\
		4	&	3	&	12	&	0.759	&	0.22	&	0.384	&	0.223	&	0.376	\\
		4	&	4	&	128	&	8.785	&	5.095	&	53.451	&	3.814	&	68.479	\\
		4	&	5	&	24	&	10.892	&	1.46	&	2.986	&	0.972	&	4.646	\\
		5	&	4	&	64	&	10.362	&	3.04	&	17.2	&	6.213	&	23.313	\\
		5	&4	&1024&	5140.657&	229.350	&3545.935&	43.481&	4331.662\\
		\hline
		%		\midrule
		
		%		\bottomrule
		%		
		
	\end{tabular}
	
\end{table}

\section{Conclusions}\label{sec:conclusions}
The premise of this paper was to develop an algorithm which would generate the minimum cardinality non-anticipativity constraint set for multistage stochastic programs where the uncertainty is endogenous and the realization of the endogenous parameter is gradual. The algorithm utilized knowledge of the subsets of scenarios formed by the realization of uncertainty in order to construct necessary and minimal sets of non-anticipativity constraints. A running time analysis showed that algorithm scaled as $O(S^3)$ as $S\rightarrow +\infty$ as long as the number of uncertain parameters is kept constant. Illustrative computational experiments confirm the scalability of the algorithm with respect to the number of scenarios and demonstrate significant reduction in the number of NACs generated and corresponding improvement in MSSP solution time.
%We also computationally investigated the performance of the algorithm. The results revealed that for the two product manufacturing example, non-anticipativity for six scenarios randomly selected from the full set of scenarios could be enforced using five NAC constraints. This result represented a 40\% reduction in the size of the NAC set. A separate set of computational studies verified that the algorithm scaled as predicted, $O(s^3)$, in the limits of $s$, and revealed that the percent reduction in size of the NAC set increased as the number of scenarios increased. 

\section*{Acknowledgments.}
We acknowledge support  by the US National Science Foundation through the Auburn University Integrative Graduate Research and Education Traineeship (IGERT) program (Award \# 1069004). 
% Enter the text of acknowledgments here

%% The Appendices part is started with the command \appendix;
%% appendix sections are then done as normal sections
\appendix

\renewcommand\nomgroup[1]{%
  \item[\bfseries
  \ifstrequal{#1}{I}{Indices and Sets}{%
  \ifstrequal{#1}{P}{Parameters}{%
  \ifstrequal{#1}{V}{Variables}{}}}%
]}

\section{ MSSP Formulation used in the study, following \citet{colvin-2008}}\label{sec:app_form}

% \noindent \textbf{\textit{Indices/Sets}}

\nomenclature[I]{$i \in I$}{drugs}

\nomenclature[I]{$j \in J$}{clinical trials; J = $\{$PI, PII, PIII$\}$}

\nomenclature[I]{$r \in R$ }{  resource type}

\nomenclature[I]{$s,sp \in S $ }{ scenarios}
\nomenclature[I]{$p \in \{1,2,\ldots T\}$ }{ time periods}

\nomenclature[I]{$t \in \{1,2,\ldots T + max{}_{ij}(\tau{}_{ij})\}$ }{ time periods}

% \nomenclature{Indices/Sets}

\nomenclature[I]{$B$ }{ pairs of scenarios that differ in the outcome of one clinical trials $(i, j)$; non-anticipativity constraints are expressed only for $(s,sp) \in B$.}

% \noindent \textbf{\textit{Parameters}}

\nomenclature[P]{${C}_{i, j}$ }{  cost of trial time periods}

\nomenclature[P]{${Cd}_{t}$ }{  discounting factor for time value of money}

\nomenclature[P]{${f}_{i, j}$ }{  discounting factor for open revenue}

\nomenclature[P]{${P}_{i, j}$ }{  probability of trial (i, j) being successful}

\nomenclature[P]{${p}_{s}$ }{  probability of scenario s}

\nomenclature[P]{${rev}_{i}^{max}$  }{ maximum possible revenue for drug i}

\nomenclature[P]{${rev}_{ij}^{open}$ }{ estimated revenue realized for drug i, if trial (i, j $-$ 1) completed while trial (i, j) not started, is successfully developed beyond the end of the time horizon}

\nomenclature[P]{${\gamma}_{i}^{D}$ }{  loss coefficient---late completion}

\nomenclature[P]{${\gamma}_{i}^{L}$  }{ loss coefficient---loss in active patent life}

\nomenclature[P]{${\rho}_{i,j,r}$  }{ resources of type r required to start trial (i, j) }

\nomenclature[P]{${\tau}_{i,\ j}$ }{  duration of trial (i, j)}

\nomenclature[P]{${\rho}_{r}^{max}$ }{  maximum resources available of type r}

% \noindent \textbf{\textit{Variables}}

\nomenclature[V]{${Cst}_{s}$  }{  total development cost in scenario s}

\nomenclature[V]{${Rv}_{s}$ }{  revenue of scenario s}

\nomenclature[V]{${FRv}_{s}$ }{  free revenue of scenario s }

\nomenclature[V]{ENPV }{  expected net present value}

 \nomenclature[V]{${X}_{i,j,p,s}$}{   1 if clinical trial $(i, j)$ starts at the time $p$ for scenario $s$}

 \nomenclature[V]{${V}_{i,j,t,s}$ }{ 1 if clinical trial $(i, j)$ is completed by the beginning of period $t$ in scenario $s$}

 \nomenclature[V]{${Z}_{i,j,t,s}$ }{ 1 if clinical trial $(i, j)$ can be started at the beginning of time $t$ for scenario $s$}

\printnomenclature

\begin{align*}
\min \quad &ENPV=\ \sum_s{p_s\left({Rv}_s+\ {FRev}_s-{Cst}_s\right)} &\\
\sto \quad & V_{i,j,t,s}=V_{i,j,t-1,s}+X_{i,j,t-{\tau }_{i,j},s},\quad\quad\quad \forall i,j,t,s\\
&Z_{i,1,1,s}=1-X_{i,1,t,s},\quad\quad\quad\forall i,s\\
&Z_{i,1,t,s}=Z_{i,1,t-1,s}-X_{i,j,t,s},\quad\quad\quad \forall i,t>1,s  \\
&Z_{i,j,t,s}=\ Z_{i,j,t-1,s}+X_{i,j-1,t-{\tau }_{i,j-1},s}-X_{i,j,t,s},\quad\quad\quad \forall i,\ j>1,t,s  \\
& \sum_t{X_{i,j,t,s}}\le 1,\quad\quad\quad \forall i,j,s \\ 
& \sum_{tt\le t}{X_{i,j,tt,s}}\le \ V_{i,j-1,t,s},\quad\quad\quad\forall i,j>1,t,s \\ 
& \sum_i{\sum_j{\sum^{t\le t}_{tt>t-{\tau }_{i,j}}{{\rho }_{i,j,r}X_{i,j,tt,s}\le {{\rho }_r}^{max}}}},\quad\quad\quad\forall r,t,s \\ 
& X_{i,1,1,s}=X_{i,1,1,1},\quad\quad\quad\forall i,s \\ 
& {-V}_{i^{s,sp},j^{s,sp},t,s\ }\le \ X_{i,j,t,s}-X_{i,j,t,sp}\le V_{i^{s,sp},j^{s,sp},t,s\ }\ \forall i,j,\ \left(s,sp\right)\in B,\quad\quad\quad t>1 \\ 
& {Cst}_s=\sum_{i,j,t}{{cd}_tC_{i,j}X_{i,j,t,s}},\quad\quad\quad\forall s \\ 
& {Rv}_s=\sum_i{\sum_t{{rev}^{max}_iX_{i,PII,t,s}-{{\gamma }_i}^D\left(Z_{i,PII,t,s}+Z_{i,PIII,t,s}\right)-{{\gamma }_i}^L(t+{\tau }_{i,PIII})X_{i,PIII,t,s}}},\quad\quad\quad\forall s \\ 
& {FRev}_s=\sum_i{\sum_j{{rev}^{open}_{i,j}f_{i,j}Z_{i,j,\left|T\right|,s}}}+\sum_i{\sum_{j\in \{PI,PII\}}{\sum_{t>\left|T\right|-{\tau }_{i,j}}{{rev}^{run}_{i,j,t}f_{i,j+1}X_{i,j,t,s}}}},\quad\quad\quad\forall s \\ 
& {{rev}_{i,j}}^{open}={{rev}_i}^{max}-{{\gamma }_i}^L\left(\left|T\right|+\sum_{jj\ge j}{{\tau }_{i,jj}}\right) \\ 
& {{rev}_{i,j}}^{run}={{rev}_i}^{max}-{{\gamma }_i}^L\left(t+\sum_{jj\ge j}{{\tau }_{i,jj}}\right) \\
& f_{i,j} = 0.9 \bigg[ \frac{rev_i^{max} - \gamma^L_i |T| - \sum_{j'\ge j} C_{i,j}}{rev_i^{max} - \gamma_i^L |T|} \bigg] \\
& X_{i,j,t,s}, V_{i,j,t,s}, Z_{i,j,t,s} \in \{ 0,1\} 
\end{align*}

\section{Case details}\label{sec:app_params}

The parameter values used in the case studies are reported in Tables \ref{t:app1}, \ref{t:app2}, \ref{t:app3}, \ref{t:app4}, \ref{t:app5}.

\begin{table}[!h!]
\caption{Parameters of the two-drug case study. Clinical trial plan for a 15-month planning horizon divided into five equal time periods.\label{t:app1}
}
\tiny	 
\begin{tabular}{|c|c|c|c|c|c|c|c|c|c|c|c|c|c|} \hline 
\textbf{Drug} & \multicolumn{2}{|p{0.4in}|}{\textbf{Duration}} & \multicolumn{2}{|p{0.9in}|}{\textbf{Probability of Success}} & \multicolumn{2}{|p{0.3in}|}{\textbf{Cost (\$M)}} & \multicolumn{2}{|p{0.5in}|}{\textbf{Resource 1 (Max = 2)}} & \multicolumn{2}{|p{0.5in}|}{\textbf{Resource 2 (Max = 3)}} & \textbf{\textit{rev${}^{max}$}} & \textbf{\textit{${\gamma}^{L}$}} & \textbf{\textit{${\gamma}^{D}$}} \\ \hline 
\textbf{} & \textbf{PI} & \textbf{PII} & \textbf{PI} & \textbf{PII} & \textbf{PI} & \textbf{PII} & \textbf{PI} & \textbf{PII} & \textbf{PI} & \textbf{PII} &  &  &  \\ \hline 
D1 & 2 & 4 & 0.3 & 0.5 & 10 & 90 & 1 & 1 & 1 & 2 & 3100 & 19.2 & 44 \\ \hline 
D2 & 2 & 3 & 0.4 & 0.6 & 10 & 80 & 1 & 2 & 1 & 1 & 3250 & 19.6 & 56 \\ \hline 
\end{tabular}

\end{table}

\begin{table}[!h!]
\tiny	
\caption{Parameters of the three-drug case study. 
	Clinical Trial plan for a 36-month planning horizon divided into 12 equal time periods.\label{t:app2}
}
\begin{tabular}{|c|c|c|c|c|c|c|c|c|c|c|c|c|c|c|c|c|c|c|} \hline 
\textbf{Drug} & \multicolumn{3}{|p{0.5in}|}{\textbf{Duration}} & \multicolumn{3}{|p{0.9in}|}{\textbf{Probability of Success}} & \multicolumn{3}{|p{0.6in}|}{\textbf{Trial Cost (\$M)}} & \multicolumn{3}{|p{0.5in}|}{\textbf{Resource 1 (Max = 2)}} & \multicolumn{3}{|p{0.5in}|}{\textbf{Resource 1 (Max = 3)}} & \textbf{\textit{rev${}^{max}$}} & \textbf{\textit{${\gamma}^{L}$}} & \textbf{\textit{${\gamma}^{D}$}} \\ \hline 
 & \textbf{PI} & \textbf{PII} & \textbf{PIII} & \textbf{PI} & \textbf{PII} & \textbf{PIII} & \textbf{PI} & \textbf{PII} & \textbf{PIII} & \textbf{PI} & \textbf{PII} & \textbf{PIII} & \textbf{PI} & \textbf{PII} & \textbf{PIII} &  &  &  \\ \hline 
D1 & 2 & 4 & 4 & 0.3 & 0.5 & 0.8 & 10 & 90 & 220 & 1 & 1 & 2 & 1 & 2 & 3 & 3100 & 19.2 & 44 \\ \hline 
D2 & 2 & 3 & 5 & 0.4 & 0.6 & 0.8 & 10 & 80 & 200 & 1 & 2 & 2 & 1 & 1 & 3 & 3250 & 19.6 & 56 \\ \hline 
D3 & 2 & 3 & 4 & 0.3 & 0.6 & 0.9 & 10 & 90 & 180 & 1 & 1 & 2 & 1 & 1 & 3 & 3300 & 20 & 52 \\ \hline 
\end{tabular}

\end{table}

\begin{table}[!h!]
	\caption{Parameters of the four-drug case study. Clinical trial plan for an 18-month planning horizon divided into six equal time periods.\label{t:app3}
	}
\tiny
\begin{tabular}{|c|c|c|c|c|c|c|c|c|c|c|c|c|c|c|c|c|c|c|} \hline 
\textbf{Drug} & \multicolumn{3}{|p{0.5in}|}{\textbf{Duration}} & \multicolumn{3}{|p{0.9in}|}{\textbf{Probability of Success}} & \multicolumn{3}{|p{0.6in}|}{\textbf{Trial Cost (\$M)}} & \multicolumn{3}{|p{0.5in}|}{\textbf{Resource 1 (Max = 4)}} & \multicolumn{3}{|p{0.5in}|}{\textbf{Resource 1 (Max = 3)}} & \textbf{\textit{rev${}^{max}$}} & \textbf{\textit{${\gamma}^{L}$}} & \textbf{\textit{${\gamma}^{D}$}} \\ \hline 
 & \textbf{PI} & \textbf{PII} & \textbf{PIII} & \textbf{PI} & \textbf{PII} & \textbf{PIII} & \textbf{PI} & \textbf{PII} & \textbf{PIII} & \textbf{PI} & \textbf{PII} & \textbf{PIII} & \textbf{PI} & \textbf{PII} & \textbf{PIII} &  &  &  \\ \hline 
D1 & 1 & 1 & 3 & 0.3 & 0.5 & 0.8 & 10 & 90 & 220 & 1 & 1 & 2 & 1 & 2 & 3 & 3100 & 19.2 & 22 \\ \hline 
D2 & 1 & 2 & 2 & 0.4 & 0.6 & 0.8 & 10 & 80 & 200 & 1 & 2 & 2 & 1 & 1 & 3 & 3250 & 19.6 & 28 \\ \hline 
D3 & 1 & 1 & 3 & 0.3 & 0.6 & 0.9 & 10 & 90 & 180 & 1 & 1 & 2 & 1 & 1 & 3 & 3300 & 20 & 26 \\ \hline 
D4 & 1 & 2 & 2 & 0.4 & 0.6 & 0.8 & 10 & 100 & 170 & 1 & 1 & 2 & 1 & 2 & 3 & 3000 & 19.4 & 24 \\ \hline 
\end{tabular}
\end{table}

\begin{table}[!h!]
	\tiny
\caption{Parameters of the five-drug case study. Clinical trial plan for an 18-month planning horizon divided into six equal time periods.\label{t:app4}
}

	\begin{tabular}{|c|c|c|c|c|c|c|c|c|c|c|c|c|c|c|c|c|c|c|} \hline \textbf{Drug} & \multicolumn{3}{|p{0.5in}|}{\textbf{Duration}} & \multicolumn{3}{|p{0.8in}|}{\textbf{Probability of Success}} & \multicolumn{3}{|p{0.6in}|}{\textbf{Trial Cost (\$M)}} & \multicolumn{3}{|p{0.5in}|}{\textbf{Resource 1 (Max = 4)}} & \multicolumn{3}{|p{0.5in}|}{\textbf{Resource 1 (Max = 3)}} & \textbf{\textit{rev${}^{max}$}} & \textbf{\textit{${\gamma}^{L}$}} & \textbf{\textit{${\gamma}^{D}$}} \\ \hline 
 & \textbf{PI} & \textbf{PII} & \textbf{PIII} & \textbf{PI} & \textbf{PII} & \textbf{PIII} & \textbf{PI} & \textbf{PII} & \textbf{PIII} & \textbf{PI} & \textbf{PII} & \textbf{PIII} & \textbf{PI} & \textbf{PII} & \textbf{PIII} &  &  &  \\ \hline 
D1 & 1 & 1 & 3 & 0.3 & 0.5 & 0.8 & 10 & 90 & 220 & 1 & 1 & 2 & 1 & 2 & 3 & 3100 & 19.2 & 22 \\ \hline 
D2 & 1 & 2 & 2 & 0.4 & 0.6 & 0.8 & 10 & 80 & 200 & 1 & 2 & 2 & 1 & 1 & 3 & 3250 & 19.6 & 28 \\ \hline 
D3 & 1 & 1 & 3 & 0.3 & 0.6 & 0.9 & 10 & 90 & 180 & 1 & 1 & 2 & 1 & 1 & 3 & 3300 & 20 & 26 \\ \hline 
D4 & 1 & 2 & 2 & 0.4 & 0.6 & 0.8 & 10 & 100 & 170 & 1 & 1 & 2 & 1 & 2 & 3 & 3000 & 19.4 & 24 \\ \hline 
D5 & 1 & 2 & 3 & 0.35 & 0.5 & 0.9 & 10 & 70 & 210 & 1 & 1 & 2 & 1 & 1 & 3 & 3150 & 19.6 & 24 \\ \hline 
\end{tabular}
\end{table}

\begin{table}[!ht!]
	\tiny
\caption{Parameters for a six-drug case study. Clinical trial plan for an 18-month planning horizon divided into six equal time periods.\label{t:app5}
}

	\begin{tabular}{|c|c|c|c|c|c|c|c|c|c|c|c|c|c|c|c|c|c|c|} \hline 
\textbf{Drug} & \multicolumn{3}{|p{0.5in}|}{\textbf{Duration}} & \multicolumn{3}{|p{0.9in}|}{\textbf{Probability of Success}} & \multicolumn{3}{|p{0.6in}|}{\textbf{Trial Cost (\$M)}} & \multicolumn{3}{|p{0.5in}|}{\textbf{Resource 1 (Max = 4)}} & \multicolumn{3}{|p{0.5in}|}{\textbf{Resource 1 (Max = 3)}} & \textbf{\textit{rev${}^{max}$}} & \textbf{\textit{${\gamma}^{L}$}} & \textbf{\textit{${\gamma}^{D}$}} \\ \hline 
 & \textbf{PI} & \textbf{PII} & \textbf{PIII} & \textbf{PI} & \textbf{PII} & \textbf{PIII} & \textbf{PI} & \textbf{PII} & \textbf{PIII} & \textbf{PI} & \textbf{PII} & \textbf{PIII} & \textbf{PI} & \textbf{PII} & \textbf{PIII} &  &  &  \\ \hline 
D1 & 1 & 1 & 3 & 0.3 & 0.5 & 0.8 & 10 & 90 & 220 & 1 & 1 & 2 & 1 & 2 & 3 & 3100 & 19.2 & 22 \\ \hline 
D2 & 1 & 2 & 2 & 0.4 & 0.6 & 0.8 & 10 & 80 & 200 & 1 & 2 & 2 & 1 & 1 & 3 & 3250 & 19.6 & 28 \\ \hline 
D3 & 1 & 1 & 3 & 0.3 & 0.6 & 0.9 & 10 & 90 & 180 & 1 & 1 & 2 & 1 & 1 & 3 & 3300 & 20 & 26 \\ \hline 
D4 & 1 & 2 & 2 & 0.4 & 0.6 & 0.8 & 10 & 100 & 170 & 1 & 1 & 2 & 1 & 2 & 3 & 3000 & 19.4 & 24 \\ \hline 
D5 & 1 & 2 & 3 & 0.35 & 0.5 & 0.9 & 10 & 70 & 210 & 1 & 1 & 2 & 1 & 1 & 3 & 3150 & 19.6 & 24 \\ \hline 
D6 & 1 & 2 & 3 & 0.45 & 0.45 & 0.8 & 10 & 85 & 195 & 1 & 2 & 2 & 2 & 1 & 3 & 3050 & 19 & 25 \\ \hline 
\end{tabular}
\end{table}

 \section*{References}
%% \label{}

%% If you have bibdatabase file and want bibtex to generate the
%% bibitems, please use
%%
  \bibliographystyle{elsarticle-harv} 
  \bibliography{NAC}

%% else use the following coding to input the bibitems directly in the
%% TeX file.

%\begin{thebibliography}{00}
%
%%% \bibitem[Author(year)]{label}
%%% Text of bibliographic item
%
%\bibitem[ ()]{}
%
%\end{thebibliography}
\end{document}